\documentclass[a4paper,10pt]{amsart}

\usepackage{amsmath,amssymb,amsthm,hyperref,graphicx,ytableau,environ,dsfont,enumitem,tikz,float,bm,stmaryrd,comment,scontents}
\newenvsc{forshort}[store-env=forshort,print-env=false]

\newtheorem{thm}{Theorem}
\newtheorem{prop}[thm]{Proposition}
\newtheorem{cor}[thm]{Corollary}
\newtheorem{lemma}[thm]{Lemma}
\newtheorem{conjecture}[thm]{Conjecture}

\theoremstyle{definition}

\newtheorem{example}[thm]{Example}
\newtheorem{definition}[thm]{Definition}
\newtheorem{remark}[thm]{Remark}

\allowdisplaybreaks

\newcommand{\ii}{\mathrm{i}}
\newcommand{\ee}{\mathrm{e}}
\newcommand{\KZ}{\mathrm{KZ}}
\newcommand{\QQ}{\mathbb{Q}}
\newcommand{\ZZ}{\mathbb{Z}}
\newcommand{\NN}{\mathbb{N}}
\newcommand{\CC}{\mathbb{C}}
\newcommand{\RR}{\mathbb{R}}
\newcommand{\Zfrak}{\mathfrak{Z}}
\newcommand{\Zfrakbar}{\bar{\mathfrak{Z}}}

\DeclareMathOperator{\diag}{diag}
\DeclareMathOperator{\SSYT}{SSYT}
\DeclareMathOperator{\sym}{\mathcal{S}}
\DeclareMathOperator{\fin}{\mathcal{A}}

\DeclareMathOperator{\symhat}{\widehat{\mathcal{S}}}
\DeclareMathOperator{\finhat}{\widehat{\mathcal{A}}}
\DeclareMathOperator{\unifhat}{\widehat{\mathcal{U}}}
\DeclareMathOperator{\Pcal}{\mathcal{P}}
\DeclareMathOperator{\sgn}{sgn}
\DeclareMathOperator{\id}{id}
\DeclareMathOperator{\type}{type}
\DeclareMathOperator{\rim}{Rim}

\DeclareMathOperator{\ned}{ned}

\DeclareMathOperator{\wt}{wt}

\begin{document}
\title{Unified Schur multiple zeta functions}
\author{Nako Kamiyama}
\address{Nako Kamiyama, Sophia University, Chiyoda City, Tokyo 102-8554, Japan}
\email{n-kamiyama-4c3@eagle.sophia.ac.jp}
\author{Simon Rutard}
\address{Simon Rutard, IMSc, CIT Campus, Taramani, Chennai 600113, India}
\email{simont@imsc.res.in}
\date{\today}
\begin{abstract}
    We define symmetrized multiple zeta values, unified zeta functions, and finite zeta values in the Schur setting. We establish a Jacobi--Trudi formula and a Giambelli formula for unified Schur multiple zeta functions, and we show that there is a one-to-one correspondence between finite and symmetrized Schur multiple zeta values under the Kaneko--Zagier conjecture. We also show a vanishing property for a symmetric sum of finite Schur multiple zeta values, generalizing a result obtained by Hoffman for finite multiple zeta values.
\end{abstract}
\maketitle
\ytableausetup{centertableaux}

\section{Introduction}

We define \emph{multiple zeta functions} (MZFs for short) by the multiple Dirichlet series
\[
    \zeta(s_{1},\ldots,s_n) := \sum_{0<m_1<\cdots<m_n} \frac{1}{m_1^{s_1} \cdots m_n^{s_n}}, \qquad (\Re s_1,\ldots,\Re s_{n-1} \geq 1, \Re s_n>1).
\]
When $s_1,\ldots,s_{n-1} \in \NN, s_n \in \ZZ_{\geq 2}$, we call $ \zeta(s_{1},\ldots,s_n)$ a \emph{multiple zeta value} (MZV for short), and we denote by $\Zfrak$ the $\QQ$-subalgebra of $\RR$ generated by all MZVs. When replacing the inequalities in the summation indices by weak inequalities, we say that such values are \emph{star MZVs}, denoted by $\zeta^{\star}(s_1,\ldots,s_n)$.

In recent years, developments have been made into defining an arithmetic analogue of MZV living in the so-called \emph{poor man's adèle ring}, which is defined by $\fin := \prod_{p \in \Pcal} (\ZZ / p\ZZ) \bigg/ \bigoplus_{p \in \Pcal} (\ZZ / p\ZZ)$ where $\Pcal$ denotes the set of prime numbers. Kaneko and Zagier conjectured (see \cite{kaneko2026finitemzv}) that there is a correspondence between \emph{finite multiple zeta values}
\[
    \zeta^{\fin}(\bm{s}):=\left( \sum_{0<m_1<\cdots<m_n<p} \frac{1}{m_1^{s_1} \cdots m_n^{s_n}} \mod p \right)_{p \in \Pcal} \in \fin \qquad (s_1,\ldots,s_n \in \ZZ)
\]
and \emph{symmetrized multiple zeta values}
\[
    \zeta^{\sym}(\bm{s}):= \sum_{k=0}^n (-1)^{s_{k+1}+\cdots+s_n} \zeta_*(s_1,\ldots,s_k)\zeta_*(s_n,\ldots,s_{k+1}) \mod (\pi^2)
\]
where $s_1,\ldots,s_n \in \NN$ and where $\zeta_*$ is the harmonic regularized MZV (see \cite{ihara2006doubleshuffle}).\footnote{It is well-known that the shuffle or the harmonic regularization yield the same equivalence class of  $\zeta^{\sym}(s_1,\ldots,s_n)$ modulo the ideal $(\pi^2)$.} One can also define in a similar fashion a star analogue $\zeta^{\sym\star}(\bm{s})$ by replacing the MZV with a star MZV. Note that symmetric values are quotient class of the space of MZVs modulo $\pi^2$, therefore we note $\Zfrakbar := \Zfrak/(\pi^2)$.

We also have a refined version of these MZVs introduced independently in \cite{jarossay2014finitesymmzv} and \cite{rosen2015finitemzv}, along with a refined Kaneko--Zagier conjecture. The formalization we consider here comes from \cite{ono2021tadicsymmmzv}. Denote by $\finhat$ the completion of $\fin$ defined by
\[
    \finhat = \varprojlim \left( \prod_{p \in \Pcal} (\ZZ/p^k\ZZ) \bigg/ \bigoplus_{p \in \Pcal} (\ZZ/p^k\ZZ) \right),
\]
where the projective maps are given by
\begin{align*}
	\prod_{p \in \Pcal} (\ZZ/p^m\ZZ) \bigg/ \bigoplus_{p \in \Pcal} (\ZZ/p^m\ZZ) &\to \prod_{p \in \Pcal} (\ZZ/p^n\ZZ) \bigg/ \bigoplus_{p \in \Pcal} (\ZZ/p^n\ZZ) \\
	(x_p \mod p^{m})_{p \in \Pcal} &\mapsto (x_q \mod p^n)_{p \in \Pcal}.
\end{align*}
for $m\geq n$. Denote by $\bm{p}=(p \mod p^k)_{p \in \Pcal,k \in \NN}$. We equip $\finhat$ with the $\bm{p}$-adic topology, which is defined as the topology induced by the projections 
\[
	\finhat \to  \prod_{p \in \Pcal} (\ZZ/p^k\ZZ) \bigg/ \bigoplus_{p \in \Pcal} (\ZZ/p^k\ZZ) \qquad (k \geq 1).
\]
We call $\finhat$\emph{-multiple zeta value} the element of $\finhat$ defined by
\[
    \zeta^{\finhat}(\bm{s}) = \left( \sum_{0<m_1<\cdots<m_n<p} \frac{1}{m_1^{s_1} \cdots m_n^{s_n}} \mod p^k \right)_{p,k} \in \finhat,
\]
and denote by $\Zfrak_{\finhat}$ the closed $\QQ$-subalgebra of $\finhat$ generated by $\finhat$-MZVs and by $\bm{p}$. We define $\symhat$\emph{-multiple zeta value} by
\begin{multline*}
    \zeta^{\symhat}(\bm{s}) := \sum_{i=0}^n (-1)^{s_{i+1}+\cdots+s_n} \zeta_*(s_1,\ldots,s_i) \\
    \times \sum_{\ell_{i+1},\ldots,\ell_n \geq 0} \left( \prod_{j=i+1}^n \binom{s_j+\ell_j-1}{\ell_j} t^{\ell_j} \right) \zeta_*(s_n+\ell_n,\ldots,s_{i+1}+\ell_{i+1}) \mod (\pi^2)
\end{multline*}
where $t$ is a formal variable, and $s_1,\ldots,s_n \in \NN$. We equip $\Zfrakbar \llbracket t \rrbracket$ with the $t$-adic topology, meaning the topology induced by the projection
\[
	\Zfrakbar \llbracket t \rrbracket \to \Zfrakbar \llbracket t \rrbracket/(t^m), \ \sum_{k \geq 0} a_k t^k \mapsto \sum_{k=0}^{m-1} a_k t^k \mod t^m \qquad (m \geq 0).
\]
The constant coefficient of $\zeta^{\symhat}(\bm{s})$ clearly corresponds to $\zeta^{\sym}(\bm{s})$. The refined Kaneko--Zagier conjecture predicts that there is a one-to-one correspondence between $\zeta^{\symhat}(\bm{s})$ and $\zeta^{\finhat}(\bm{s})$ for $s_1,\ldots,s_n \in \NN$.

\begin{conjecture}[Refined Kaneko--Zagier Conjecture]
	There exists a continuous $\QQ$-algebra isomorphism $\varphi_{\KZ} : \Zfrak_{\finhat} \to \Zfrakbar\llbracket t \rrbracket$ such that $\varphi_{\KZ}(\bm{p})=t$, and $\varphi_{\KZ}(\zeta^{\finhat}(\bm{s})) = \zeta^{\symhat}(\bm{s})$ for all $s_1,\ldots,s_n \in \NN$.
\end{conjecture}

A recent preprint by Shin-Ichiro Seki seems to establish a conditional proof of this conjecture \cite{seki2026refinedkz} under the dimension conjecture and the direct sum conjecture for ordinary and finite MZVs.

Komori \cite{komori2018finitemzv} provides another point of view to these refinements by introducing an analytic function called a \emph{unified multiple zeta function} defined by
\begin{equation}
    \zeta^{\unifhat}(s_1,\ldots,s_n;t) := \sum_{k=0}^n \ee^{\pi \ii(s_{k+1}+\cdots+s_n)}\zeta(s_{1},\ldots,s_k)\zeta(s_n,\ldots,s_{k+1};t), \label{eq:def_umzf}
\end{equation}
where we set
\[
    \zeta(s_{1},\ldots,s_n;t) := \sum_{0<m_1<\cdots<m_n} \frac{1}{(m_1-t)^{s_1} \cdots (m_n-t)^{s_n}} \qquad (t \in \CC \setminus \NN).
\]
We consider all throughout the text that the logarithm is defined on the principal branch, and we write $(-z)^{s} = \ee^{-\pi \ii s} \ee^{s \ln z}$ for all $\Re z>0$ and $s \in \CC$.

We define a star analogue $\zeta^{\unifhat\star}$ after replacing $\zeta(s_1,\ldots,s_n)$ by $\zeta^{\star}(s_1,\ldots,s_n)$, and $\zeta(s_1,\ldots,s_n;t)$ by $\zeta^{\star}(s_1,\ldots,s_n;t) = \sum_{0<m_1 \leq \cdots \leq m_n} \frac{1}{(m_1-t)^{s_1} \cdots (m_n-t)^{s_n}}$ in~\eqref{eq:def_umzf} respectively.

When $\bm{s}$ is an integer tuple, Komori showed that we have an expansion near $t=0$ of the form
\[
    \zeta^{\unifhat}(\bm{s};t) = \sum_{k \geq 0} \zeta^{\unifhat,k}(\bm{s}) t^k,
\]
with $\zeta^{\unifhat,k}(\bm{s}) \in \Zfrak[\pi \ii]$.

Note that there is an isomorphism of $\QQ$-algebras
\[
    \Zfrak[\pi \ii]\llbracket t \rrbracket/(\pi \ii) \xrightarrow{\sim} \Zfrak\llbracket t \rrbracket/(\pi^2), \ \sum_{k=0}^{\infty} z_k t^k \mapsto \sum_{k=0}^{\infty} \Re z_k t^k,
\]
and for the rest of the article, we identify both algebras. At positive integers $s_1,\ldots,s_n \in \NN$, Komori relates the special values $\zeta^{\unifhat}(\bm{s};t)$ with $\symhat$-MZV. More specifically, when $t$ is treated as a formal variable, he established that the class of $\zeta^{\unifhat}(\bm{s};t)$ modulo $\pi \ii \Zfrak[\pi \ii] \llbracket t \rrbracket$ is equal to the class of $\zeta^{\symhat}(\bm{s})$ modulo $\pi^2\Zfrak \llbracket t \rrbracket$. Therefore, one can consider $\zeta^{\unifhat}(\bm{s};t)$ as an analytic analogue of $\symhat$-MZV. For the rest of the paper, we denote $\Zfrakbar := \Zfrak[\pi \ii] / \pi \ii \Zfrak[\pi \ii]$.

Komori showed that UMZF are holomorphic for $(\bm{s},t) \in \CC^{n} \times \CC\setminus \NN$. For $s_1,\ldots,s_n \in \ZZ$, we set
$\zeta^{\symhat}(\bm{s}) := \zeta^{\unifhat}(\bm{s};t) \mod \pi \ii \Zfrak[\pi \ii] \llbracket t \rrbracket$, and we set $\zeta^{\symhat\star}(\bm{s}) := \zeta^{\unifhat \star}(\bm{s};t) \mod \pi \ii \Zfrak[\pi \ii] \llbracket t \rrbracket$.

\begin{thm}[\cite{ono2023refinedkanekozagier}]
\label{th:unified_hurwitz}
    Assuming the refined Kaneko--Zagier conjecture, there exists a continuous $\QQ$-algebra isomorphism $\varphi_{\KZ} \colon \Zfrak_{\finhat} \xrightarrow{\sim} \Zfrakbar\llbracket t \rrbracket$ such that $\varphi_{\KZ}(\zeta^{\finhat}(\bm{s})) = \zeta^{\symhat}(\bm{s})$ holds for all indices $\bm{s} \in \ZZ^r$.
\end{thm}

\begin{remark}
	Note that one can also define symmetrized multiple zeta values by considering limits of truncated sums (see \cite{zhao2016mzvmpl}). We denote by $\prec$ the \emph{strict Kontsevich order} which is a strict total order on $\ZZ_{\neq 0}$ defined by $1 \prec 2 \prec \cdots \prec ``+\infty \prec -\infty" \prec \cdots \prec -2 \prec -1$, and we denote by $\preceq$ the weak Kontsevich order (meaning that $i \preceq j$ if and only if $i \prec j$ or $i=j$). We have
	\[
	    \zeta^{\sym}(s_1,\ldots,s_n)=\lim_{N \to +\infty} \sum_{\substack{m_1, \ldots, m_n \in \llbracket -N,N\rrbracket\setminus \{0\} \\ m_1 \prec \cdots \prec m_n}} \frac{1}{m_1^{s_1} \cdots m_n^{s_n}} \mod \pi^2.
	\]
\end{remark}

Komori showed that we can retrieve multiple harmonic sums, and obtained
\[
    \zeta^{\fin}(k_1,\ldots,k_n)=\left(\zeta^{\unifhat}(-k_1^{(p)},\ldots,-k_n^{(p)};p) \mod p\right)_{p \in \Pcal} \in \fin \qquad (k_1,\ldots,k_n \in \ZZ)
\]
where $\zeta^{\unifhat}(-k_1^{(p)},\ldots,-k_n^{(p)};p):=\lim_{\varepsilon \to 0} \zeta^{\unifhat}(-k_1^{(p)},\ldots,-k_n^{(p)};p-\varepsilon)$, and
\[
    k^{(p)} := 
    \begin{cases}
        p-1-k & \text{ if } k>0, \\
        -k & \text{ else}.
    \end{cases}
\]

Let 
\begin{equation}
    \zeta^{\pm N}(\bm{s};t):= \sum_{\substack{m_1, \ldots, m_n \in \llbracket-N,N\rrbracket\setminus \{0\} \\ m_1 \prec \cdots \prec m_n}} \frac{1}{(m_1+\mathds{1}_{<0}(m_1)t)^{s_1} \cdots (m_n+\mathds{1}_{<0}(m_n)t)^{s_n}}. \label{eq:trun_zeta_Kontsevich}
\end{equation}
We can split the above sum into pieces of the form $\{ (m_1,\ldots,m_n) \mid 0<m_1<\cdots<m_j, m_{j+1}<\cdots<m_n<0\}$ in order to have $\zeta^{\unifhat}(\bm{s};t)=\lim_{N \to +\infty} \zeta^{\pm N}(\bm{s};t)$ for $\Re s_1, \ldots, \Re s_n >1$.

Similarly, using the weak Kontsevich order we can define a star version of the unified multiple zeta functions by first setting
\begin{equation}
    \zeta^{\pm N \star}(\bm{s};t):= \sum_{\substack{m_1, \ldots, m_n \in \llbracket-N,N\rrbracket\setminus \{0\} \\ m_1 \preceq \cdots \preceq m_n}} \frac{1}{(m_1+\mathds{1}_{<0}(m_1)t)^{s_1} \cdots (m_n+\mathds{1}_{<0}(m_n)t)^{s_n}}. \label{eq:trun_zetastar_Kontsevich}
\end{equation}
Then it is easy to show that $\zeta^{\unifhat\star}(\bm{s};t)=\lim_{N \to +\infty} \zeta^{\pm N\star}(\bm{s};t)$ for $\Re s_1,\ldots, \Re s_n > 1$.

Throughout the paper, we denote by $\NN$ the set of positive integers, by $\ZZ$ the set of integers, and we write $\llbracket n,m \rrbracket := \{ n,n+1,\ldots,m\}$. We denote by $\mathfrak{S}_n$ the symmetric group of $n$ elements. The empty sum should be understood as $0$, the empty product as $1$, and $F(\emptyset)$ as $1$ for a function $F$.

\section{Finite, symmetrized and unified analogues for Schur multiple zeta functions}

We first recall the definition of \emph{finite Schur multiple zeta values}. Then, we introduce a truncated Schur multiple zeta function (with respect to the Kontsevich order), generalizing~\eqref{eq:trun_zeta_Kontsevich} and~\eqref{eq:trun_zetastar_Kontsevich}. By taking its limit, we define a unified Schur multiple zeta function. Throughout the following sections, we fix a partition $\lambda$ of a fixed integer $n$, and we consider $T(\lambda,X)$ to be the set of all Young tableaux of shape $\lambda$ over a set $X$. We denote by $D(\lambda)$ the Young diagram of $\lambda$. We associate to a partition $\lambda$ its transpose $\lambda'$.

In the following, we denote by 
\[
	\zeta_{\lambda}^{N}(\bm{s}) = \sum_{M \in \SSYT(\lambda,\llbracket 1,N \rrbracket)} \frac{1}{M^{\bm{s}}}
\]
the truncated Schur multiple zeta function, as defined in \cite{nakasuji2018schurgeneralization}. We introduce finite and unified analogues for Schur multiple zeta functions.

\begin{definition}
   Let $\bm{s} \in T(\lambda,\ZZ)$. We define $\finhat$\emph{-Schur multiple zeta values} by
    \[
        \zeta^{\finhat}_{\lambda}(\bm{s}):=\left( \zeta^{p-1}_{\lambda}(\bm{s}) \mod p^k \right)_{p \in \Pcal, k} \in \finhat.
    \]
    and we define \emph{finite Schur multiple zeta values} by
    \[
        \zeta^{\fin}_{\lambda}(\bm{s}):=\left( \zeta^{p-1}_{\lambda}(\bm{s}) \mod p \right)_{p \in \Pcal} \in \fin.
    \]
\end{definition}

Note that for all prime numbers $p$, the denominator of the truncation $\zeta^{p-1}_{\lambda}(\bm{s})$ is prime to $p$ because the finite sum features denominators that are prime to $p$, therefore finite and $\finhat$-Schur multiple zeta values are well-defined.

\begin{definition}
    Let $N \in \ZZ_{>0}$ be an integer. We denote by $\SSYT^{\pm N}(\lambda)$ the finite set of all tuples $M=(m_{ij})_{(i,j) \in D(\lambda)} \in T(\lambda,\{-N,\ldots,-1,1,\ldots,N\})$ such that we have $m_{i1} \preceq m_{i2} \preceq \cdots$ for all $i$, and $m_{1j} \prec m_{2j} \prec \cdots$ for all $j$. For any $t \in \CC\setminus \NN$, we denote $M(t):=(m_{ij}+\mathds{1}_{<0}(m_{ij})t)_{(i,j) \in D(\lambda)}$ where $\mathds{1}_{<0}$ denotes the indicator function of the set $\RR_{<0}$. We now define the \emph{truncated unified Schur multiple zeta function} by considering the finite sum
    \[
        \zeta_{\lambda}^{\pm N}(\bm{s};t) := \sum_{M \in \SSYT^{\pm N}(\lambda)} \frac{1}{M(t)^{\bm{s}}}
    \]
    where $\bm{s}=(s_{ij})_{(i,j) \in D(\lambda)}$ is a tuple of complex numbers, and $M(t)^{\bm{s}} = \prod_{(i,j) \in D(\lambda)} (m_{ij}+\mathds{1}_{<0}(m_{ij})t)^{s_{ij}}$. We denote the \emph{weight} of $\bm{s}$ by $\wt(\bm{s}) := \sum_{(i,j) \in D(\lambda)} s_{ij}$. We define (when it exists) the \emph{unified Schur multiple zeta function} by the limit
    \[
    	\zeta_{\lambda}^{\unifhat}(\bm{s};t) = \lim_{N \to \infty} \zeta_{\lambda}^{\pm N}(\bm{s};t).
    \]
\end{definition}

\begin{prop}
    We set
    \[
    W_{\lambda} := \{\bm{s}=(s_{ij})\in T(\lambda,\CC) \mid \Re s_{ij}>1 \text{ for all } (i,j) \in D(\lambda)\}.
    \]
    Then, the following limit exist for all $\bm{s}\in W_\lambda$ and all $t \in \CC \setminus \NN$,
    \[
    	\sum_{M \in \SSYT^{\pm N}(\lambda)} \left|\frac{1}{M(t)^{\bm{s}}}\right|.
    \]
\end{prop}

\begin{proof}
    We assume $\Re s_{ij}>1$ for all $(i,j) \in D(\lambda)$. We clearly have the following inequality
    \[
        |(m_{ij}-t)^{-s_{ij}}| \leq (C(t) m_{ij})^{-\Re s_{ij}} \ee^{2\pi |\Im s_{ij}|},
    \]
    where $C(t) := \inf_{m \geq 1} \tfrac{|m-t|}{m}>0$. By summing over all $M \in \SSYT^{\pm N}(\lambda)$, we get
    \begin{align*}
        \sum_{M \in \SSYT^{\pm N}(\lambda)} \left|\frac{1}{M(t)^{\bm{s}}} \right| &\leq \sum_{M \in \SSYT^{\pm N}(\lambda)} \prod_{(i,j) \in D(\lambda)} \ee^{2\pi |\Im s_{ij}|} |C(t)m_{ij}|^{-\Re s_{ij}} \\
        &\leq \prod_{(i,j) \in D(\lambda)} \left(1+C(t)^{-\Re s_{ij}} \ee^{2\pi |\Im s_{ij}|}\right) \zeta(\Re s_{ij})<\infty.
    \end{align*}
\end{proof}

\begin{example}
    We denote by $\lambda=(1^n)=(1,\ldots,1)$ (resp. $\lambda=(n)$) the one-column (resp. one-row) partition of $n$. Then for $\bm{s}=(s_1,\ldots,s_n)$ we have
    \[
        \zeta^{\unifhat}_{(1^n)}(\bm{s};t)=\zeta^{\unifhat}(\bm{s};t), \qquad \zeta^{\unifhat}_{(n)}(\bm{s};t)=\zeta^{\unifhat\star}(\bm{s};t).
    \]
\end{example}

In the following theorem, we show that unified Schur multiple zeta functions are holomorphic on $(\bm{s},t) \in \CC^{|\lambda|} \times \CC \setminus \NN$. Therefore, for $\bm{s}$ an integer tuple, we define \emph{symmetrized Schur multiple zeta values} by
\[
\zeta^{\sym}_{\lambda}(\bm{s}) := \zeta^{\unifhat}_{\lambda}(\bm{s};0) \mod \pi \ii \Zfrak[\pi \ii]
\]
and we define $\symhat$\emph{-Schur multiple zeta values}\footnote{It is not to be mistaken with symmetric Schur multiple zeta functions, coined by Maki Nakasuji and Wataru Takeda in an unrelated paper.} by
\[
\zeta^{\symhat}_{\lambda}(\bm{s}) := \zeta^{\unifhat}_{\lambda}(\bm{s};t) \mod \pi \ii \Zfrak[\pi \ii] \llbracket t \rrbracket.
\]

\begin{thm}
	Let $\lambda$ be a partition of $n$. The unified Schur multiple zeta function $\zeta_{\lambda}^{\unifhat}(\bm{s};t)$ is holomorphic for $(\bm{s},t) \in \CC^{n} \times \CC \setminus \NN$. Furthermore, at integer tuples we have $\zeta_{\lambda}^{\unifhat}(\bm{s};t) \in \Zfrak[\pi \ii] \llbracket t \rrbracket$ and $\zeta_{\lambda}^{\finhat}(\bm{s}) \in \Zfrak_{\finhat}$, and assuming the Kaneko--Zagier conjecture, there is a correspondence between $\symhat$-Schur MZV and $\finhat$-Schur MZV,
	\[
	\forall \bm{s} \in \ZZ^{n}, \qquad \varphi_{\KZ}(\zeta^{\finhat}_{\lambda}(\bm{s})) = \zeta^{\symhat}_{\lambda}(\bm{s}).
	\]
\end{thm}

\begin{proof}
	The holomorphy follows from the fact that a USMZF can also be written as a finite sum of UMZFs, just like a Schur MZF can be written down in terms of MZFs as established in \cite{nakasuji2018schurgeneralization}. 
	
	Let $\mathcal{F}(\lambda)$ be the set of all bijections $f:D(\lambda)\to \llbracket 1,n \rrbracket$ satisfying the following two conditions:
	\begin{enumerate}[label=(\roman*)]
		\item for all $i$, $f((i,j))<f((i,j'))$ if and only if $j<j'$, 
		\item for all $j$, $f((i,j))<f((i',j))$ if and only if $i<i'$.
	\end{enumerate} 
	Moreover, for $\bm{s}=(s_{ij})\in T(\lambda,\CC)$, put
	\[
	V(\bm{s})=
	\left\{ \left(s_{f^{-1}(1)},\ldots,s_{f^{-1}(n)}\right)\in \CC^{n} \mid f\in \mathcal{F}(\lambda)
	\right\}.
	\] 
	We write $\bm{u} \sqsubseteq \bm{s}$ for $\bm{u}=(u_1,\ldots,u_m)\in \CC^m$ if there exists $(v_1,\ldots,v_{n})\in V(\bm{s})$ satisfying the following: for all $1\le k\le m$, there exist $1\le h_k\le m$ and $l_k\ge 0$ such that  
	\begin{enumerate}[label=(\roman*)]
		\item $u_k=v_{h_k}+\cdots +v_{h_k+l_k}$,
		\item there are no $i$ and $i'$ such that $i\ne i'$ and $s_{ij},s_{i'j}\in\{v_{h_k},\ldots ,v_{h_k+l_k}\}$ for some $j$,
		\item $\bigsqcup^{m}_{k=1}\{h_k,\ldots,h_k+l_k\}=\llbracket 1,n \rrbracket$.
	\end{enumerate}
	Then, by the definition, we have
	\[
	\zeta^{\pm N}_{\lambda}(\bm{s};t) = \sum_{\bm{u} \sqsubseteq \bm{s}}\zeta^{\pm N}(\bm{u};t).
	\]
	Moreover, by an inclusion-exclusion argument as in \cite{nakasuji2018schurgeneralization}, one can also obtain the ``dual" expressions
	\[
	\zeta^{\pm N}_{\lambda}(\bm{s};t) = \sum_{\bm{u} \sqsubseteq \bm{s}'}(-1)^{n-\ell(\bm{u})}\zeta^{\pm N \star}(\bm{u};t)
	\]
	where $\bm{s}' \in T(\lambda',W_{\lambda'})$ is such that $s'_{ij} = s_{ji}$. By taking the limit $N \to +\infty$ for $\bm{s} \in T(\lambda,W_{\lambda})$, we obtain
	\begin{align*}
		\zeta^{\unifhat}_{\lambda}(\bm{s};t) =& \sum_{\bm{u} \sqsubseteq \bm{s}}\zeta^{\unifhat}(\bm{u};t), \\
		=& \sum_{\bm{u} \sqsubseteq \bm{s}'}(-1)^{n-\ell(\bm{u})}\zeta^{\unifhat \star}(\bm{u};t),
	\end{align*}
	where we write $\ell(\bm{u}) := m$ when $\bm{u}=(u_1,\ldots,u_m)$. We know from \cite{komori2021unified} that $\zeta^{\unifhat}(u_1,\ldots,u_m;t)$ is entire on $\CC^{m} \times \CC \setminus \NN$. Each term $\bm{u} \sqsubseteq \bm{s}$ can be expressed as components in $\bm{s}=(s_{ij})$, therefore unified Schur multiple zeta functions are holomorphic functions on $(\bm{s},t) \in \CC^{|\lambda|} \times \CC \setminus \NN$. It follows from Theorem~\ref{th:unified_hurwitz} that $\zeta_{\lambda}^{\symhat}(\bm{s}) \in \Zfrak$ for $\bm{s} \in \ZZ^{|\lambda|}$.
	
	Let $\bm{s} \in \ZZ^{|\lambda|}$. In \cite{nakasuji2018schurgeneralization}, they show
	\[
	\zeta^N_{\lambda}(\bm{s}) = \sum_{\bm{u} \sqsubseteq \bm{s}} \zeta^N(\bm{u}).
	\]
	By considering the sequence modulo $p^k$ for $p$ a prime and $k \in \mathbb{N}$, we get an equality between $\finhat$-MZVs and $\finhat$-Schur MZVs
	\[
	\zeta^{\finhat}_{\lambda}(\bm{s}) = \sum_{\bm{u} \sqsubseteq \bm{s}} \zeta^{\finhat}(\bm{u}).
	\]
	In particular, $\zeta^{\finhat}_{\lambda}(\bm{s}) \in \Zfrak_{\finhat}$.
	
	Assuming the refined Kaneko--Zagier conjecture, from Theorem~\ref{th:unified_hurwitz}, we obtain
	\begin{align*}
		\varphi_{\KZ}(\zeta^{\finhat}_{\lambda}(\bm{s})) &= \sum_{\bm{u} \sqsubseteq \bm{s}} \varphi_{\KZ}(\zeta^{\finhat}(\bm{u})) \\
		&= \sum_{\bm{u} \sqsubseteq \bm{s}} \zeta^{\symhat}(\bm{u}) \\
		&= \zeta^{\symhat}_{\lambda}(\bm{s}).
	\end{align*}
\end{proof}

\section{Main results}

We study some analytic properties of USMZFs, and we obtain that the unified version satisfy a determinant formula. We already know from \cite{nakasuji2018schurgeneralization} that truncated Schur multiple zeta values satisfy a determinant formula, and it implies that $\finhat$-Schur multiple zeta values also satisfy a determinant formula.

We now establish several determinant formulas that USMZF and $\finhat$-Schur MZV satisfy. First, we denote the conjugate of a partition $\lambda=(\lambda_1,\ldots,\lambda_r)$ by $\lambda'=(\lambda'_1,\ldots,\lambda'_s)$. The determinant formulas will be satisfied for tuples $\bm{s}$ belonging to $T^{\diag}(\lambda,\CC)=\{\bm{s}=(s_{ij})\in T(\lambda,\CC) \mid s_{ij}=s_{kl} \text{ if } j-i=l-k\}$. For $\bm{s}=(s_{ij})\in T^{\diag}(\lambda,\CC)$, we always write $a_k=s_{i,i+k}$. For example, when $\lambda=(3,3,1)$,
$\bm{s}=(s_{ij})\in T^{\diag}((3,3,1),\CC)$ implies that $\bm{s}$ is of the form
\[
\bm{s}=
\ytableausetup{boxsize=18pt,aligntableaux=center}
\begin{ytableau}
	s_{11} & s_{12} & s_{13} \\
	s_{21} & s_{22} & s_{23} \\
	s_{31}
\end{ytableau} 
=
\begin{ytableau}
	a_0 & a_1 & a_2 \\
	a_{-1} & a_0 & a_1 \\
	a_{-2}
\end{ytableau} 
\, .
\]\\
In both determinant formulas, we assume that the value of an MZF at a tuple is $1$ if the length of the tuple (i.e. the number of components of the tuple) is $0$, and $0$ if the length of the tuple is negative.
\begin{thm}
    \label{th:jacobi_trudi_formula}
    Let $\lambda=(\lambda_1,\ldots,\lambda_r)$ be a partition, and let $(\bm{s},t) \in T^{\diag}(\lambda,\CC) \times \CC \setminus \NN$. Then $\zeta_{\lambda}^{\unifhat}(\bm{s};t)$ satisfies two determinant formulas,
	\begin{align}
		\zeta_{\lambda}^{\unifhat}(\bm{s};t) &= \det\left[\zeta^{\unifhat\star}(a_{-j+1},a_{-j+2},\ldots,a_{-j+(\lambda_{i}-i+j)};t)\right]_{1\le i,j\le r}, \label{eq:jacobi_trudi_formula_Htype} \\
		\zeta_{\lambda}^{\unifhat}(\bm{s};t) &= \det\left[\zeta^{\unifhat}(a_{j-1},a_{j-2},\ldots,a_{j-(\lambda'_{i}-i+j)};t)\right]_{1\le i,j\le s}. \label{eq:jacobi_trudi_formula_Etype}
	\end{align}
\end{thm}

Since USMZF are holomorphic for $(\bm{s},t) \in \CC^{|\lambda|} \times \CC \setminus \NN$, it is enough to prove the previous theorem on $\bm{s} \in W_{\lambda}^{\diag}$ and $(\Re t<1)$, where
\[
	W^{\diag}_{\lambda}
	=W_\lambda\cap T^{\diag}(\lambda,\CC)
\]
and then extend the result to $(\bm{s},t) \in T^{\diag}(\lambda,\CC) \times \CC \setminus \NN$ by analytic continuation.

We now write the two determinant formulas satisfied by $\finhat$-Schur MZVs. First, we write
\[
	\zeta^{\finhat \star}(s_1,\ldots,s_n) := \left(\sum_{0<m_1\leq\cdots\leq m_n} \frac{1}{m_1^{s_1} \cdots m_n^{s_n}} \mod p^k\right)_{p \in \Pcal, k \in \NN}.
\]

\begin{thm}
    Let $\lambda=(\lambda_1,\ldots,\lambda_r)$ be a partition, and let $\bm{s} \in T^{\diag}(\lambda,\ZZ)$. Then $\zeta^{\fin}_{\lambda}(\bm{s})$ satisfies two determinant formulas,
    \begin{align*}
        \zeta_{\lambda}^{\finhat}(\bm{s}) &= \det\left[\zeta^{\finhat \star}(a_{-j+1},a_{-j+2},\ldots,a_{-j+(\lambda_{i}-i+j)})\right]_{1\le i,j\le r}, \\
        \zeta_{\lambda}^{\finhat}(\bm{s}) &= \det\left[\zeta^{\finhat}(a_{j-1},a_{j-2},\ldots,a_{j-(\lambda'_{i}-i+j)})\right]_{1\le i,j\le s}.
    \end{align*}
\end{thm}

This last theorem follows from the Jacobi--Trudi formula for truncated Schur multiple zeta functions established in \cite{nakasuji2018schurgeneralization}, along with the next lemma.

\begin{lemma}
	Let $(R_k, \varphi_{k})_{k \in \NN}$ be a projective system of commutative rings, and let $R = \varprojlim R_k$. Let $n \in \ZZ_{\geq 0}$, and let $a = (a(k))_{k \in \NN} \in R$ and $b_{r,s} = (b_{i,j}(k))_{k \in \NN} \in R$ ($1 \leq i,j \leq n$). Suppose there exists $k_0 \in \NN$ such that for all $k \geq k_0$, we have
	\[
	a(k) = \det (b_{i,j}(k))_{1 \leq i,j \leq n} \in R_k.
	\]
	Then $a = \det (b_{i,j})_{1 \leq i,j \leq n}$
\end{lemma}

\begin{proof}
    By definition of the determinant we have
    \begin{align*}
        \det (b_{i,j})_{1 \leq i,j \leq n} &= \sum_{\sigma \in \mathfrak{S}_n} \sgn(\sigma) \prod_{i=1}^n b_{i,\sigma(i)} \\
        & = \left( \sum_{\sigma \in \mathfrak{S}_n} \sgn(\sigma) \prod_{i=1}^n b_{i,\sigma(i)}(k) \right)_{k \in \NN} \\
        &= (\det(b_{i,j}(k))_{1 \leq i,j \leq n})_{k \in \NN} \\
        &= (a(k))_{k \in \NN}.
    \end{align*}
\end{proof}

Let $\lambda$ be a partition, and denote by $\mathfrak{S}_{\lambda}$ the permutation group of the cells of the Young diagram $D(\lambda)$. We consider its left action on $T(\lambda,\CC)$, denoted by $\sigma \cdot \bm{s}$ ($\sigma \in \mathfrak{S}_{\lambda}, \bm{s} \in T(\lambda,\ZZ)$), defined by $\sigma \cdot s_{ij} = s_{\sigma^{-1}(ij)}$. The next theorem shows that the symmetric sum of finite Schur MZVs vanishes.

\begin{thm}
    \label{TH:SYMMETRIC_RELATION_FINITE_SCHUR_MZV}
     Let $\lambda$ be a partition of $n \geq 1$. Then for all $\bm{s} \in T(\lambda,\NN)$, we have
    \begin{equation}
        \sum_{\sigma \in \mathfrak{S}_{\lambda}} \zeta^{\fin}_{\lambda}(\sigma \cdot \bm{s}) = 0. \label{eq:symmetric_relation_finite_schur_mzv}
    \end{equation}
\end{thm}

Under the Kaneko--Zagier conjecture, we would also have
\[
	\sum_{\sigma \in \mathfrak{S}_{\lambda}} \zeta^{\sym}_{\lambda}(\sigma \cdot \bm{s}) = 0 \qquad (\bm{s} \in T(\lambda,\ZZ_{\geq 2})).
\]

\begin{remark}
    When $\lambda$ is a one-row or one-column Young tableau, the previous result is known as Hoffman's formula \cite[Theorem~4.4]{hoffman2015quasisymm}. In the case where $\lambda$ is a hook shape partition, then one can use the determinant formula satisfied for finite Schur MZVs together with Hoffman's formula to show~\eqref{eq:symmetric_relation_finite_schur_mzv}.
\end{remark}

\begin{remark}
	This equality does not seem to be true if we replace finite Schur MZV with $\finhat$-Schur MZV. Based on SageMath computations, an interesting phenomenon appears when computing the reduction modulo $p^k$ (for a fixed $k \geq 2$) of the symmetric sum,
	\[
		\sum_{\sigma \in \mathfrak{S}_{\lambda}} \zeta^{p-1}_{\lambda}(\sigma \cdot \bm{s}) \mod p^k.
	\]
	We experimentally found that there seem to be a vanishing for the symmetric sum modulo $p^k$ where the primes $p$ are big enough, for $k$ smaller than the size of the biggest square in the Young tableau.
\end{remark}

\section{Proof of Theorem~\ref{TH:SYMMETRIC_RELATION_FINITE_SCHUR_MZV}}

We denote by $\mu(M)=(\mu_1,\ldots,\mu_N)$ the \emph{content} of $M \in \SSYT(\lambda,\llbracket 1,N \rrbracket)$, meaning that $i$ appears $\mu_i$ times in $M$, for every $i=1,\ldots,N$.

Let $P_k(x_1,\ldots,x_N) = x_1^{k} + \cdots + x_N^{k}$ denote the $k$-th symmetric power sum in $N$ variables, and let $m_{\mu} = \sum_{\sigma \in \mathfrak{S}_N} x_{\sigma(1)}^{\mu_1} \cdots x_{\sigma(N)}^{\mu_N}$ denote the \emph{monomial symmetric polynomial}. Recall both families form a basis for the vector space of symmetric polynomials in $N$-variables. The \emph{Kostka number} $K_{\lambda, \mu}$ denotes the number of Young tableaux $M \in \SSYT(\lambda,\llbracket 1,N \rrbracket)$ with content equal to $\mu$. They satisfy
\begin{equation}
	s_{\lambda}(x_1,\ldots,x_N) = \sum_{\mu} K_{\lambda,\mu} m_{\mu}(x_1,\ldots,x_N). \label{eq:schur_decomposed_into_monomial}
\end{equation}

Let
\[
    Q^N_{\lambda,\bm{s}}(x_1,\ldots,x_N) := \sum_{M \in \SSYT(\lambda,\llbracket 1,N \rrbracket)} \prod_{c \in D(\lambda)} x^{s_c}_{M(c)}
\]
and consider its \emph{symmetrization}
\begin{align*}
    F^N_{\lambda,\bm{s}}(x_1,\ldots,x_N) &:= \sum_{\sigma \in \mathfrak{S}_{\lambda}} Q^N_{\lambda,\sigma \cdot \bm{s}}(x_1,\ldots,x_N) \\
    &= \sum_{\sigma \in \mathfrak{S}_{\lambda}} \sum_{M \in \SSYT(\lambda,\llbracket 1,N \rrbracket)} \prod_{c \in D(\lambda)} x^{s_{\sigma(c)}}_{M(c)}.
\end{align*}
First observe that truncated Schur multiple zeta values (at $N=p-1$) corresponds to $Q^{p-1}_{\lambda,\bm{s}}(1,\ldots,\tfrac{1}{p-1})$, and that $F^{p-1}_{\lambda,\bm{s}}(1,\ldots,\tfrac{1}{p-1})$ corresponds to the symmetric sum at the left-hand side of~\eqref{eq:symmetric_relation_finite_schur_mzv}. We first show that the polynomial $F^N_{\lambda,\bm{s}}(x_1,\ldots,x_N)$ is symmetric in $x_1,\ldots,x_N$.

\begin{lemma}
    The polynomial $F^N_{\lambda,\bm{s}}(x_1,\ldots,x_N)$ is a homogeneous symmetric polynomial of degree $\wt(\bm{s})$.
\end{lemma}

\begin{proof}    
    Let $M,M' \in \SSYT(\lambda,\llbracket 1,N \rrbracket)$ such that they have the same content. Then there is a permutation $\tau \in \mathfrak{S}_{\lambda}$ such that $M(c) = M'(\tau(c))$ for all $c \in D(\lambda)$. We have that $x^{s_c}_{M'(c)} = x^{s_c}_{M(\tau^{-1}(c))}$, and by doing a change of variable of the form $\sigma'=\sigma \tau^{-1}$, we get
    \begin{align*}
        \sum_{\sigma \in \mathfrak{S}_{\lambda}} \prod_{c \in D(\lambda)} x^{s_{\sigma(c)}}_{M(c)} &= \sum_{\sigma \in \mathfrak{S}_{\lambda}} \prod_{c \in D(\lambda)} x^{s_{\sigma(\tau(c))}}_{M(c)} \\
        &= \sum_{\sigma \in \mathfrak{S}_{\lambda}} \prod_{c \in D(\lambda)} x^{s_{\sigma(c)}}_{M(\tau^{-1}(c))} \\
        &= \sum_{\sigma \in \mathfrak{S}_{\lambda}} \prod_{c \in D(\lambda)} x^{s_{\sigma(c)}}_{M'(c)}.
    \end{align*}
    For a fixed $\mu$, and for $M$ with content $\mu$, we set
    \[
        H^N_{\lambda, \mu, \bm{s}}(x_1,\ldots,x_N) := \sum_{\sigma \in \mathfrak{S}_{\lambda}} \prod_{c \in D(\lambda)} x^{s_{\sigma(c)}}_{M(c)}.
    \]
    Grouping the sum over $M \in \SSYT(\lambda,\llbracket 1,N \rrbracket)$ according to content, we can write
    \begin{align*}
    	F^N_{\lambda,\bm{s}}(x_{\tau(1)},\ldots,x_{\tau(N)}) &= \sum_{\mu} K_{\lambda,\mu}
    	H^N_{\lambda,\tau(\mu),\bm{s}}(x_1,\ldots,x_N) \\
    	&= \sum_{\nu} K_{\lambda,\nu} H^N_{\lambda,\nu,\bm{s}}(x_1,\ldots,x_N) \\
    	&= F^N_{\lambda,\bm{s}}(x_1,\ldots,x_N),
    \end{align*}
    Let $\tau \in \mathfrak{S}_N$, and let $\mu=(\mu_1,\ldots,\mu_N)$. We denote $\tau(\mu) := (\mu_{\tau^{-1} 1},\ldots,\mu_{\tau^{-1} N})$, and we observe that if $M$ has content $\mu$, then $\tau M$ has content $\tau(\mu)$. By symmetry of Schur polynomials and by \eqref{eq:schur_decomposed_into_monomial}, it is clear that $K_{\lambda, \mu} = K_{\lambda, \tau(\mu)}$. Also, by a change of variable we obtain
    \begin{align*}
        \sum_{\mu(M) = \mu} \sum_{\sigma \in \mathfrak{S}_{\lambda}} \prod_{c \in D(\lambda)} x^{s_{\sigma(c)}}_{\tau(M(c))} = \sum_{\mu(M) = \tau(\mu)} \sum_{\sigma \in \mathfrak{S}_{\lambda}} \prod_{c \in D(\lambda)} x^{s_{\sigma(c)}}_{M(c)},
    \end{align*}
    which shows that $F^N_{\lambda,\bm{s}}(\tau(x_1,\ldots,x_N)) = F^N_{\lambda,\bm{s}}(x_1,\ldots,x_N)$.
\end{proof}

Since $F_{\lambda,\bm{s}}^N$ is a symmetric homogeneous polynomial, we can write it down in terms of a polynomial (without constant terms) in the symmetric power sums. We set $N=p-1$ where $p$ is a prime. Now it remains to see that when $x_1=1,\ldots,x_{p-1} = \tfrac{1}{p-1}$, we find
\[
    P_k\left(1,\ldots,\frac{1}{p-1}\right) = \frac{1}{1^k} + \cdots + \frac{1}{(p-1)^k}.
\]
It is well-known that the last sum vanishes modulo $p$ (when $p>k+1$), as a consequence of the Seki--Bernoulli formula (see for example \cite[\S1]{bachmann2026mzv}).

\section{Determinant formula}

The proof follows essentially from \cite{nakasuji2018schurgeneralization}, where we change the patterns and the weights. Instead of considering paths inside the lattice $\ZZ^2$, with vertical coordinates going from $1$ to $N \in \NN$, we consider paths with vertical coordinates going from $1$ to $N$, then from $-N$ to $-1$. We also consider weights translated by the parameter $t$, in order to match our definition of unified Schur multiple zeta function. We first prove a determinant formula for the truncated sum $\zeta^{\pm N}_{\lambda}(\bm{s};t)$, and then we take the limit as $N \to +\infty$.

In §\ref{subsection:rim_decomp}, we recall the vocabulary introduced in \cite{nakasuji2018schurgeneralization}. In §\ref{subsection:patterns}, we define the patterns in a lattice suited for the Kontsevich order. We also define weights shifted by a parameter $t$ to fit this new setting.

\subsection{Rim decompositions of a partition}
\label{subsection:rim_decomp}

A skew partition is a pair of partitions $(\lambda,\mu)$ satisfying $\mu\subset \lambda$, that is $\mu_i\le \lambda_i$ for all $i$. The resulting skew shape is denoted by $\lambda/\mu$ and the corresponding Young diagram is denoted by $D(\lambda/\mu)$. We often identify $\lambda/\mu$ with $D(\lambda/\mu)$. A skew Young diagram $\theta$ is called a \emph{ribbon} if $\theta$ is connected and contains no $2\times 2$ block of boxes.

Let $\lambda$ be a partition. A sequence $\Theta=(\theta_1,\ldots,\theta_t)$ of ribbons is called a \emph{rim decomposition of $\lambda$} if $D(\theta_k)\subset D(\lambda)$ for $1\le k\le t$ and $\theta_1\sqcup\cdots\sqcup\theta_k$, the gluing of $\theta_1,\ldots,\theta_k$, is (the Young diagram of) a partition $\lambda^{(k)}$ for $1\le k\le t$ satisfying $\lambda^{(t)}=\lambda$.

One can naturally identify a rim decomposition $\Theta=(\theta_1,\ldots,\theta_t)$ of $\lambda$ with the Young tableau $T=(t_{ij})\in T(\lambda,\llbracket 1,t \rrbracket)$ defined by $t_{ij}=k$ if $(i,j)\in D(\theta_k)$. 
\begin{example}
    \label{ex:rim}
    The following $\Theta=(\theta_1,\theta_2,\theta_3,\theta_4)$ is a rim decomposition of $\lambda=(4,3,3,2)$
    \[
        \ytableausetup{boxsize=normal,aligntableaux=center}
        \Theta
        =
        \begin{ytableau}
            1 & 1 & 3 & 3 \\
            2 & 3 & 3 \\
            2 & 3 & 4 \\
            3 & 3 
        \end{ytableau}
    \]
    which means that
    \ytableausetup{mathmode,boxsize=10pt,aligntableaux=center}
        $\theta_1=\ydiagram{2}$\,,
        $\theta_2=\ydiagram{0,1,1}$\,,
        $\theta_3=\ydiagram{2+2,1+2,1+1,2}$\,, and
        $\theta_4=\ydiagram{1}$\,.
\end{example}

Write $\lambda=(\lambda_1,\ldots,\lambda_r)$. We call a rim decomposition $\Theta=(\theta_1,\ldots,\theta_r)$ of $\lambda$ an \emph{$H$-rim decomposition} if each $\theta_i$ starts at $(i,1)$ for all $1\le i\le r$. Here, we permit $\theta_i=\emptyset$. We denote by $\rim^{\lambda}_H$ the set of all $H$-rim decompositions of $\lambda$.

\begin{example}
    \label{ex:Hrim}
    The following $\Theta=(\theta_1,\theta_2,\theta_3,\theta_4)$ is an $H$-rim decomposition of $\lambda=(4,3,3,2)$
    \[
        \ytableausetup{boxsize=normal,aligntableaux=center}
        \Theta=
        \begin{ytableau}
            1 & 1 & 3 & 3 \\
            3 & 3 & 3 \\
            3 & 4 & 4 \\
            4 & 4 
        \end{ytableau}
    \]
    which means that
    \ytableausetup{mathmode,boxsize=10pt,aligntableaux=center} 
    $\theta_1=\ydiagram{2}$\,,
    $\theta_2=\emptyset$\,, 
    $\theta_3=\ydiagram{2+2,3,1}$\, and 
    $\theta_4=\ydiagram{1+2,2}$\,.
    Note that the rim decomposition appearing in Example~\ref{ex:rim} is not an $H$-rim decomposition.
\end{example}

\subsection{Patterns}
\label{subsection:patterns}

Fix a positive integer $N$. We set 
\[
    \Lambda^{+N}:= \{ (i,j) \in \ZZ^2 \mid 1 \leq j \leq N \}, \qquad \Lambda^{-N}:= \{ (i,j) \in \ZZ^2 \mid -N \leq j \leq 0 \}.
\]
We now consider the union $\Lambda^{\pm N} := \Lambda^{+N} \sqcup \Lambda^{-N}$. We call the \emph{height} of an integer $j \in \{ -N,\ldots,N \}$ the number
\[
    h_N(j) :=
    \begin{cases}
        j &\text{if } 1 \leq j \leq N, \\
        2N+1+j &\text{if } -N \leq j \leq 0.
    \end{cases}
\]
On $\Lambda^{\pm N}$ we define the following partial order
\[
    (i_1,j_1) \preceq (i_2,j_2) \Leftrightarrow i_1 \leq i_2 \ \text{ and } h_N(j_1) \leq h_N(j_2).
\]
We call the ordered set $\Lambda^{\pm N}$ the \emph{Kontsevich horizontal strip of height $N$}.

For a partition $\lambda=(\lambda_1,\ldots,\lambda_r)$, let $A_i$ and $B_i$ be points in $\Lambda^{\pm N}$ respectively given by $A_i=(r+1-i,1)$ and $B_i=(r+1-i+\lambda_i,-1)$ for $1\le i\le r$. Put $A=(A_1,\ldots,A_r)$ and $B=(B_1,\ldots,B_r)$. An \emph{$H$-pattern} corresponding to $\lambda$ is a tuple $L=(l_1,\ldots,l_r)$ of directed paths on $\Lambda^{\pm N}$, whose directions are allowed only to go one to the right or one up, such that $l_i$ starts at $A_i$ and ends at $B_{\sigma(i)}$ for some $\sigma\in \mathfrak{S}_r$. We call such $\sigma\in \mathfrak{S}_r$ the \emph{type} of $L$ and denote it by $\sigma=\type(L)$. Note that the type of an $H$-pattern does not depend on $N$ in the sense that the number of horizontal edges of each directed path of the $H$-pattern is independent of $N$. The number of horizontal edges appearing in the path $l_i$ is called the \emph{horizontal distance} of $l_i$ and is denoted by $\rim(l_i)$. When $\type(L)=\sigma$, we simply write $L:A\to B^{\sigma}$ where $B^{\sigma}=(B_{\sigma(1)},\ldots,B_{\sigma(r)})$ and $l_i:A_i\to B_{\sigma(i)}$. It is easy to see that $\rim(l_i)=\lambda_{\sigma(i)}-\sigma(i)+i$ and $\sum^{r}_{i=1}\rim(l_i)=|\lambda|$.
 
Let $\mathcal{H}^{\pm N}_{\lambda}$ be the set of all $H$-patterns corresponding to $\lambda$ and $\mathfrak{S}^{\lambda}_H = \{\type(L)\in \mathfrak{S}_r\mid L\in \mathcal{H}^{\pm N}_{\lambda}\}$. The following is a key lemma of our study.
 
\begin{lemma}
    For $\Theta=(\theta_1,\ldots,\theta_r) \in \rim^{\lambda}_H$, there exists $L=(l_1,\ldots,l_r)\in \mathcal{H}^{\pm N}_{\lambda}$ such that $\rim(l_i)=|\theta_i|$ for all $1\le i\le r$. Moreover, the map $\tau_H:\rim^{\lambda}_H\to \mathfrak{S}^{\lambda}_H$ given by $\tau_H(\Theta)=\type(L)$ is a bijection.
\end{lemma}

\begin{proof}
	The proof comes from relabelling the rows from $1,\ldots,N,-N,\ldots,-1$ to $1,2,\ldots,2N$, which preserves the starting and the ending point of the paths, and then to use Lemma 3.4 of \cite{nakasuji2018schurgeneralization}.
\end{proof}

\begin{example}
    Let $\lambda=(4,3,3,2)$. Then, we have $\tau_H(\Theta)=(1243)\in \mathfrak{S}_4$ where $\Theta$ is the $H$-rim decomposition of $\lambda$ appearing in Example~\ref{ex:Hrim}.
\end{example}

\subsubsection{Weight of patterns}

Fix $\bm{s}=(s_{ij})\in T(\lambda,\CC)$. We assign a weight to $L=(l_1,\ldots,l_r) \in \mathcal{H}^{\pm N}_{\lambda}$ via an $H$-rim decomposition of $\lambda$ as follows. Take $\Theta=(\theta_1,\ldots,\theta_r)\in\rim^{\lambda}_{H}$ such that $\tau_H(\Theta)=\type(L)$. Then, when the $k$th horizontal edge of $l_i$ is on the $j$th row, we weight it with $\frac{1}{j^{s_{pq}}}$ if $j>0$, and with $\frac{1}{(j+t)^{s_{pq}}}$ if $j<0$, where $(p,q)\in D(\lambda)$ is the $k$th component of $\theta_i$. Now, the weight $w^{\pm N}_{\bm{s};t}(l_i)$ of the path $l_i$ is defined to be the product of weights of all horizontal edges along $l_i$. Here, we understand that $w^{\pm N}_{\bm{s};t}(l_i)=1$ if $\theta_i=\emptyset$. Moreover, we define the weight $w^{\pm N}_{\bm{s};t}(L)$ of $L\in \mathcal{H}^{\pm N}_{\lambda}$ by
\[
    w^{\pm N}_{\bm{s};t}(L)=\prod^{r}_{i=1}w^{\pm N}_{\bm{s};t}(l_i).
\]

\begin{example}
    Let $\lambda=(4,3,3,2)$. Consider $L=(l_1,l_2,l_3,l_4)\in \mathcal{H}^{\pm 2}_{(4,3,3,2)}$.
    \begin{figure}[ht]
        \begin{center}
            \begin{tikzpicture} 
                \node at (0.5,1) {$1$};
                \node at (0.5,2) {$2$};
                \node at (0.5,3) {$-2$};
                \node at (0.5,4) {$-1$};
                \node at (1,0) {$A_4$};  
                \node at (2,0) {$A_3$};
                \node at (3,0) {$A_2$};
                \node at (4,0) {$A_1$};
                \node at (1,0.5) {$(1,1)$};  
                \node at (2,0.5) {$(2,1)$};
                \node at (3,0.5) {$(3,1)$};
                \node at (4,0.5) {$(4,1)$};
                \node at (3,4.5) {$(3,-1)$};
                \node at (5,4.5) {$(5,-1)$};
                \node at (6,4.5) {$(6,-1)$}; 
                \node at (8,4.5) {$(8,-1)$};
                \node at (3,5) {$B_4$};
                \node at (5,5) {$B_3$};
                \node at (6,5) {$B_2$};
                \node at (8,5) {$B_1$};
                \node at (1,1) {$\bullet$};
                \node at (1,2) {$\bullet$};
                \node at (1,3) {$\bullet$};
                \node at (1,4) {$\bullet$};
                \node at (2,1) {$\bullet$};
                \node at (2,2) {$\bullet$};
                \node at (2,3) {$\bullet$};  
                \node at (2,4) {$\bullet$};
                \node at (3,1) {$\bullet$};
                \node at (3,2) {$\bullet$};
                \node at (3,3) {$\bullet$};  
                \node at (3,4) {$\bullet$};
                \node at (4,1) {$\bullet$};
                \node at (4,2) {$\bullet$};
                \node at (4,3) {$\bullet$};  
                \node at (4,4) {$\bullet$};
                \node at (5,1) {$\bullet$};
                \node at (5,2) {$\bullet$};
                \node at (5,3) {$\bullet$};  
                \node at (5,4) {$\bullet$};
                \node at (6,1) {$\bullet$};
                \node at (6,2) {$\bullet$};
                \node at (6,3) {$\bullet$}; 
                \node at (6,4) {$\bullet$};
                \node at (7,1) {$\bullet$};
                \node at (7,2) {$\bullet$};
                \node at (7,3) {$\bullet$}; 
                \node at (7,4) {$\bullet$};
                \node at (8,1) {$\bullet$};
                \node at (8,2) {$\bullet$};
                \node at (8,3) {$\bullet$};  
                \node at (8,4) {$\bullet$};
                \draw (4,1) -- (5,1) -- (5,2) -- (6,2) -- (6,3) -- (6,4);
                \draw[dotted] (3,1) -- (3,2) -- (3,3) -- (3,4);
                \draw[loosely dashdotted] (2,1) -- (3,1) -- (4,1) -- (4,2) -- (4,3) -- (5,3) -- (6,3) -- (7,3) -- (7,4) -- (8,4);      
                \draw[dashed] (1,1) -- (1,2) -- (2,2) -- (2,3) -- (3,3) -- (4,3) -- (4,4) -- (5,4);
            \end{tikzpicture}
        \end{center}
        \caption{$L=(l_1,l_2,l_3,l_4)\in \mathcal{H}^{\pm 2}_{(4,3,3,2)}$}
    \end{figure}
    
    Since $\type(L)=(1243)$, the corresponding $H$-rim decomposition of $\lambda$ is the one appearing in Example~\ref{ex:Hrim}. Let
    \(
        \ytableausetup{boxsize=18pt,aligntableaux=center}
        \bm{s}=
        \begin{ytableau}
            a & b & c & d \\
            e & f & g \\
            h & i & j \\
            k & l 
        \end{ytableau}
        \in T((4,3,3,2),\CC)
    \). 
    Then, the weights of the $l_i$ are given by
    \begin{align*}
        &w^{\pm 2}_{\bm{s};t}(l_1)
        =\frac{1}{1^a2^b}, \quad
        w^{\pm 2}_{\bm{s};t}(l_2)
        =1, \quad
        w^{\pm 2}_{\bm{s};t}(l_3)
        =\frac{1}{1^h1^e(-2+t)^f(-2+t)^g(-2+t)^c(-1+t)^d}, \\
        &w^{\pm 2}_{\bm{s};t}(l_4)
        =\frac{1}{2^k(-2+t)^l(-2+t)^i(-1+t)^j}.
    \end{align*}
    In particular, when 
    \(
        \ytableausetup{boxsize=18pt,aligntableaux=center}
        \bm{s}=
        \begin{ytableau}
        a_0 & a_1 & a_2 & a_3 \\
        a_{-1} & a_{0} & a_1 \\
        a_{-2} & a_{-1} & a_0 \\
        a_{-3} & a_{-2}
        \end{ytableau}
        \in T^{\diag}((4,3,3,2),\CC)
    \),
    these are equal to 
    \begin{align*}
        &w^{\pm 2}_{\bm{s};t}(l_1)
        =\frac{1}{1^{a_0}2^{a_1}}, \quad 
        w^{\pm 2}_{\bm{s};t}(l_2)
        =1, \\
        &w^{\pm 2}_{\bm{s};t}(l_3)
        =\frac{1}{1^{a_{-2}}1^{a_{-1}}(-2+t)^{a_0}(-2+t)^{a_1}(-2+t)^{a_2}(-1+t)^{a_3}}, \\
        &w^{\pm 2}_{\bm{s};t}(l_4)
        =\frac{1}{2^{a_{-3}}(-2+t)^{a_{-2}}(-2+t)^{a_{-1}}(-1+t)^{a_0}}.
    \end{align*}
    Notice that, in this case, from the definition of the weight, the tuple of indices of the exponent of the denominator of $w^{\pm 2}_{\bm{s};t}(l_i)$ along $l_i$ is $(a_{1-i},a_{1-i+1},a_{1-i+2},\ldots)$ for all $i$. 
\end{example}

\subsubsection{Proof}

A proof of Jacobi--Trudi formula~\eqref{eq:jacobi_trudi_formula_Etype} is given by calculating the sum 
\[
    X^{\pm N}_{\lambda}(\bm{s};t)
    =\sum_{L\in \mathcal{H}^{\pm N}_{\lambda}}\varepsilon_{\type(L)}w^{\pm N}_{\bm{s};t}(L)
    =\sum_{\sigma\in \mathfrak{S}^{\lambda}_H}
    \varepsilon_{\sigma}
    \sum_{L:A\to B^{\sigma}}
    w^{\pm N}_{\bm{s};t}(L),
\]
where $\varepsilon_{\sigma}$ is the signature of $\sigma\in \mathfrak{S}_r$. First, the inner sum can be calculated as follows.

\begin{lemma} 
\label{lemma:H_step_1}
    For $\sigma\in \mathfrak{S}^{\lambda}_H$, let $\Theta^{\sigma}=(\theta^{\sigma}_1,\ldots,\theta^{\sigma}_r)\in \rim^{\lambda}_H$ be the $H$-rim decomposition such that $\tau_H(\Theta^{\sigma})=\sigma$. Then, we have  
    \[ 
        \sum_{L:A\to B^{\sigma}}w^{\pm N}_{\bm{s};t}(L)
        =\prod^{r}_{i=1}\zeta^{\pm N\star}\left(\theta^{\sigma}_i(\bm{s});t\right).
    \]
    Here, for $\Theta=(\theta_1,\ldots,\theta_r)\in\rim^{\lambda}_H$, $\theta_i(\bm{s})\in\CC^{|\theta_i|}$ is the tuple obtained by reading the entries of the restriction of $\bm{s}$ to $\theta_i$ from the bottom left to the top right.
\end{lemma}

\begin{proof}
    Let $L=(l_1,\ldots,l_r)\in \mathcal{H}^{\pm N}_{\lambda}$ be an $H$-pattern of type $\sigma$. Then $l_i$ is a path from $A_i$ to $B_{\sigma(i)}$ with $\rim(l_i)=\lambda_{\sigma(i)}-\sigma(i)+i=|\theta^{\sigma}_i|$. For simplicity, write $k_i=\lambda_{\sigma(i)}-\sigma(i)+i$ and $\theta^{\sigma}_i(\bm{s})=(s_{i,1},\ldots,s_{i,k_i})$. Suppose that $l_i$ has $n_j$ steps on the $j$th row for $j \in \llbracket-N,N\rrbracket \setminus \{0\}$. Then, from the definition of the weight, we have
    \begin{align*}
        w^{\pm N}_{\bm{s};t}(l_i)
        =& \underbrace{\frac{1}{1^{s_{i,1}}}\cdots \frac{1}{1^{s_{i,n_1}}}}_{n_1 \text{ terms}}
        \cdots 
        \underbrace{\frac{1}{N^{s_{i,n_1+\cdots+n_{N-1}+1}}}\cdots \frac{1}{N^{s_{i,n_1+\cdots+n_{N}}}}}_{n_N \text{ terms}} \\
        & \cdot \underbrace{\frac{1}{(-N+t)^{s_{i,n_1+\cdots+n_{N}+1}}}\cdots \frac{1}{(-N+t)^{s_{i,n_1+\cdots+n_{N}+n_{-N}}}}}_{n_{-N} \text{ terms}}
        \\ & \cdots \underbrace{\frac{1}{(-1+t)^{s_{i,n_1+\cdots+n_{-2}+1}}}\cdots \frac{1}{(-1+t)^{s_{i,n_1+\cdots+n_{-1}}}}}_{n_{-1} \text{ terms}}
    \end{align*}
    with $n_1+\cdots+n_N+n_{-N}+\cdots+n_{-1}=k_i$. This shows that
    \begin{align*}
        \sum_{L:A\to B^{\sigma}}w^{\pm N}_{\bm{s};t}(L)
        &=\prod^{r}_{i=1}\sum_{l_i:A_i\to B_{\sigma(i)}}w^{\pm N}_{\bm{s};t}(l_i) \\
        &=\prod^{r}_{i=1}\sum_{\substack{-N \leq m_1, \ldots, m_{k_i} \leq N \\ m_1 \preceq \cdots \preceq m_{k_i} \\ m_1,\ldots,m_{k_i} \neq 0}} \frac{1}{(m_1+\mathds{1}_{<0}(m_1)t)^{s_{i,1}}\cdots (m_{k_i}+\mathds{1}_{<0}(m_{k_i})t)^{s_{i,k_i}}} \\
        &=\prod^{r}_{i=1}\zeta^{\pm N\star}(s_{i,1},\ldots,s_{i,k_i};t).
    \end{align*}
\end{proof}

From Lemma~\ref{lemma:H_step_1}, we have
\begin{equation}
\label{eq:X_N}
    X^{\pm N}_{\lambda}(\bm{s};t)
    =\sum_{\sigma\in \mathfrak{S}^{\lambda}_H}
    \varepsilon_{\sigma}
    \prod^{r}_{i=1}\zeta^{\pm N\star}\left(\theta^{\sigma}_i(\bm{s});t\right).
\end{equation}
Let $\mathcal{H}^{\pm N}_{\lambda,0}$ be the set of all $L=(l_1,\ldots,l_r)\in \mathcal{H}^{\pm N}_{\lambda}$ such that no distinct paths $l_i$ and $l_j$ intersect. Define 
\[ 
    X^{\pm N}_{\lambda,0}(\bm{s};t)=\sum_{L\in \mathcal{H}^{\pm N}_{\lambda,0}} \varepsilon_{\type(L)} w^{\pm N}_{\bm{s};t}(L),
    \quad X^{\pm N}_{\lambda,1}(\bm{s};t) =\sum_{L\in \mathcal{H}^{\pm N}_{\lambda}\setminus \mathcal{H}^{\pm N}_{\lambda,0}}\varepsilon_{\type(L)}w^{\pm N}_{\bm{s};t}(L).
\]
Clearly we have $X^{\pm N}_{\lambda}(\bm{s};t)=X^{\pm N}_{\lambda,0}(\bm{s};t)+X^{\pm N}_{\lambda,1}(\bm{s};t)$. Moreover, $\type(L)=\id$ for all $L\in \mathcal{H}^{\pm N}_{\lambda,0}$ where $\id$ is the identity element of $\mathfrak{S}_r$ and $\id$ corresponds to the trivial $H$-rim decomposition $(\theta_1,\ldots,\theta_r)=((\lambda_1),\ldots,(\lambda_r))$. Employing the well-known identification between non-intersecting lattice paths and semi-standard Young tableaux, we have
\[
    X^{\pm N}_{\lambda,0}(\bm{s};t)=\sum_{L\in \mathcal{H}^{\pm N}_{\lambda,0}}w^{\pm N}_{\bm{s};t}(L)
    =\zeta^{\pm N}_{\lambda}(\bm{s};t).
\]
Therefore, from~\eqref{eq:X_N}, we reach the expression 
\begin{equation}
\label{eq:X_main_expression}
    \zeta^{\pm N}_{\lambda}(\bm{s};t)=\sum_{\sigma\in \mathfrak{S}^{\lambda}_H} \varepsilon_{\sigma} \prod^{r}_{i=1} \zeta^{\pm N\star}\left(\theta^{\sigma}_i(\bm{s});t\right)-X^{\pm N}_{\lambda,1}(\bm{s};t).
\end{equation}

Now, one can obtain the determinant formula by taking the limit $N \to \infty$ of~\eqref{eq:X_main_expression} under suitable assumptions on $\bm{s}$ described in Theorem~\ref{th:jacobi_trudi_formula} together with the following proposition.

\begin{prop}
\label{prop:jacobi_trudi_formula_X}
    Assume that $\bm{s}=(s_{ij})\in T^{\diag}(\lambda,\CC)$.
    Write $a_{k}=s_{i,i+k}$ for $k\in\ZZ$.
    \begin{enumerate}
        \item
        We have 
        \begin{equation}
            \label{eq:jacobi_trudi_formula_X}
            X^{\pm N}_{\lambda}(\bm{s};t)
            =\det
            \left[
            \zeta^{\pm N\star}(a_{-j+1},a_{-j+2},\ldots,a_{-j+(\lambda_i-i+j)};t)
            \right]_{1\le i,j\le r}.
        \end{equation}
        \item
        It holds that
        \begin{equation}
            \label{eq:X_1_vanish}
            X^{\pm N}_{\lambda,1}(\bm{s};t)=0.
        \end{equation}
    \end{enumerate}
\end{prop}

\begin{proof}
    We first notice that, if $\bm{s}=(s_{ij})\in T^{\diag}(\lambda,\CC)$, then we have
    \[
        \theta^{\sigma}_i(\bm{s})
        =(a_{1-i},a_{1-i+1},\ldots,a_{1-i+(\lambda_{\sigma(i)}-\sigma(i)+i)-1})
    \]
    for all $1\le i\le r$. Therefore, understanding that $\zeta^{\pm N}_{(k)}=0$ for $k<0$, from~\eqref{eq:X_N}, we have
    \begin{align*}
        X^{\pm N}_{\lambda}(\bm{s};t)
        &= \sum_{\sigma\in \mathfrak{S}_{H}^{\lambda}} \varepsilon_{\sigma}
        \prod^{r}_{i=1}\zeta^{\pm N\star}\left(\theta^{\sigma}_i(\bm{s});t\right)\\
        &=\sum_{\sigma\in \mathfrak{S}_r}
        \varepsilon_{\sigma}
        \prod^{r}_{i=1}
        \zeta^{\pm N\star}(a_{1-i},a_{1-i+1},\ldots,a_{1-i+(\lambda_{\sigma(i)}-\sigma(i)+i)-1};t)\\
        &=\det
        \left[
        \zeta^{\pm N\star}(a_{1-i},a_{1-i+1},\ldots,a_{1-i+(\lambda_j-j+i)-1};t)
        \right]_{1\le i,j\le r}\\
        &=\det
        \left[
        \zeta^{\pm N\star}(a_{-j+1},a_{-j+2},\ldots,a_{-j+(\lambda_i-i+j)};t)
        \right]_{1\le i,j\le r}.
    \end{align*}
    Hence, we obtain~\eqref{eq:jacobi_trudi_formula_X}. 
    
    We next show the second assertion. To do that, we employ the well-known involution $L\mapsto \overline{L}$ on $\mathcal{H}^{\pm N}_{\lambda}\setminus \mathcal{H}^{\pm N}_{\lambda,0}$ defined as follows. For $L=(l_1,\ldots,l_r)\in \mathcal{H}^{\pm N}_{\lambda}\setminus \mathcal{H}^{\pm N}_{\lambda,0}$ of type $\sigma$, consider the first (rightmost) intersection point appearing in $L$, at which two paths, say $l_i$ and $l_j$ cross. Let $\overline{L}$ be an $H$-pattern that contains every path in $L$ except for $l_i$ and $l_j$ and two more paths $\overline{l_i}$ and $\overline{l_j}$. Here, $\overline{l_i}$ (resp. $\overline{l_j}$) follows $l_i$ (resp. $l_j$) until it meets the first intersection point and after that follows $l_j$ (resp. $l_i$) to the end. Notice that, if $\bm{s}=(s_{ij})\in T^{\diag}(\lambda,\CC)$, then we have $w^{\pm N}_{\bm{s};t}(\overline{L})=w^{\pm N}_{\bm{s};t}(L)$ since there is no change of horizontal edges between $L$ and $\overline{L}$. Moreover, we have $\type(\overline{L})=\tau \type(L)$ (where $\tau$ denotes a transposition that corresponds to the endpoints of $L$ and $\overline{L}$) because the end points of $L$ and $\overline{L}$ are just switched. These imply that 
    \begin{align*}
        X^{\pm N}_{\lambda,1}(\bm{s};t)
        &=\sum_{L\in \mathcal{H}^{\pm N}_{\lambda}\setminus \mathcal{H}^{\pm N}_{\lambda,0}}\varepsilon_{\type(\overline{L})}w^{\pm N}_{\bm{s};t}(\overline{L})\\
        &=-\sum_{L\in \mathcal{H}^{\pm N}_{\lambda}\setminus \mathcal{H}^{\pm N}_{\lambda,0}}\varepsilon_{\type(L)}w^{\pm N}_{\bm{s};t}(L)\\
        &=-X^{\pm N}_{\lambda,1}(\bm{s};t),
    \end{align*} 
    and therefore leads to~\eqref{eq:X_1_vanish}.
\end{proof}

\begin{remark}
    When $\bm{s}\in T^{\diag}(\lambda,\CC)$,~\eqref{eq:X_main_expression} can also be written in terms of the $H$-rim decomposition as follows:
    \begin{equation}
        \label{eq:main_expression_Hrim}
        \zeta^{\pm N}_{\lambda}(\bm{s};t)
        =\sum_{\Theta=(\theta_1,\theta_2,\ldots,\theta_r)\in \rim^{\lambda}_H}
        \varepsilon_H(\Theta)
        \zeta^{\pm N\star}\left(\theta_1(\bm{s});t\right)\cdots \zeta^{\pm N\star}\left(\theta_r(\bm{s});t\right),
    \end{equation}
    where $\varepsilon_H(\Theta)=\varepsilon_{\tau_H(\Theta)}$. Note that $\varepsilon_H(\Theta)=(-1)^{n-\#\{i \mid \theta_i\ne \emptyset\}}$ when $\lambda=(1^n)$.
\end{remark}

\subsection{A proof of the determinant formula of \texorpdfstring{$E$}{E}-type}

To prove~\eqref{eq:jacobi_trudi_formula_Htype}, we need to consider another type of pattern on the $\Lambda^{\pm N}$ strip. Because the discussion is essentially the same as in the previous subsection, we omit all proofs of the results in this subsection.

Let $\lambda=(\lambda_1,\ldots,\lambda_r)$ be a partition and denote by $\lambda'=(\lambda'_1,\ldots,\lambda'_s)$ its conjugate. A rim decomposition $\Theta=(\theta_1,\ldots,\theta_s)$ of $\lambda$ is called an \emph{$E$-rim decomposition} if each $\theta_i$ starts at $(1,i)$ for all $1 \le i \le s$.  Here, we again permit $\theta_i=\emptyset$. We denote by $\rim^{\lambda}_E$ the set of all $E$-rim decompositions of $\lambda$.

\begin{example}
\label{ex:E_rim}
    The following $\Theta=(\theta_1,\theta_2,\theta_3,\theta_4)$ is an $E$-rim decomposition of $\lambda=(4,3,3,2)$
    \[
        \ytableausetup{boxsize=normal,aligntableaux=center}
        \Theta =
        \begin{ytableau}
        1 & 2 & 3 & 4 \\
        1 & 2 & 3 \\
        1 & 3 & 3 \\
        3 & 3 
        \end{ytableau}
    \]
    which means that
    \ytableausetup{mathmode,boxsize=10pt,aligntableaux=center} 
    $\theta_1=\ydiagram{1,1,1}$\,,
    $\theta_2=\ydiagram{1,1}$\,, 
    $\theta_3=\ydiagram{2+1,2+1,1+2,2}$\, and 
    $\theta_4=\ydiagram{1}$\,.
\end{example}

Fix $N\in\NN$. Let $C_i$ and $D_i$ be points in $\Lambda^{\pm N}$ respectively given by $C_i=(s+1-i,1)$ and $D_i=(s+1-i+\lambda'_i,0)$ for $1\le i\le s$. Put $C=(C_1,\ldots,C_s)$ and $D=(D_1,\ldots,D_s)$. An \emph{$E$-pattern} corresponding to $\lambda$ is a tuple $L=(l_1,\ldots,l_s)$ of directed paths on $\Lambda^{\pm N}$, whose directions are allowed only to go one to the northeast or one up, such that $l_i$ starts at $C_i$ and ends at $D_{\sigma(i)}$ for some $\sigma\in \mathfrak{S}_s$. We also call such $\sigma\in \mathfrak{S}_s$ the \emph{type} of $L$ and denote it by $\sigma=\type(L)$. The number of northeast edges appearing in the path $l_i$ is called the \emph{northeast distance} of $l_i$ and is denoted by $\ned(l_i)$. When $\type(L)=\sigma$, we simply write $L:C\to D^{\sigma}$ where $D^{\sigma}=(D_{\sigma(1)},\ldots,D_{\sigma(s)})$ and $l_i:C_i\to D_{\sigma(i)}$.  It is easy to see that $\ned(l_i)=\lambda'_{\sigma(i)}-\sigma(i)+i$ and $\sum^{s}_{i=1}\ned(l_i)=|\lambda|$.

Let $\mathcal{E}^{\pm N}_{\lambda}$ be the set of all $E$-patterns corresponding to $\lambda$ and $\mathfrak{S}^{\lambda}_E=\{\type(L)\in \mathfrak{S}_s \mid L\in \mathcal{E}^{\pm N}_{\lambda}\}$. Notice that this set can be empty if $N$ is small compared to $|\lambda|$. To avoid this eventuality, we assume $2N \geq |\lambda|$ for the rest of the subsection.

\begin{lemma}
\label{lemma:rim_type_E}
    For $\Theta=(\theta_1,\ldots,\theta_s)\in \rim^{\lambda}_E$, there exists $L=(l_1,\ldots,l_s)\in \mathcal{E}^{\pm N}_{\lambda}$  such that $\ned(l_i)=|\theta_i|$ for all $1\le i\le s$. Moreover, the map $\tau_E:\rim^{\lambda}_E\to \mathfrak{S}^{\lambda}_E$ given by $\tau_E(\Theta)=\type(L)$ is a bijection.
\end{lemma}

Fix $\bm{s}=(s_{ij})\in T(\lambda,\CC)$. A weight on $L=(l_1,\ldots,l_s) \in \mathcal{E}^{\pm N}_{\lambda}$ is similarly defined via the $E$-rim decomposition of $\lambda$ as follows. Take $\Theta=(\theta_1,\ldots,\theta_s)\in\rim^{\lambda}_{E}$ such that $\tau_E(\Theta)=\type(L)$. Then, when the $k$th northeast edge of $l_i$ runs from the $j$th row to the $(j+1)$th row, we weight it with $\frac{1}{j^{s_{pq}}}$ when $j>0$, and $\frac{1}{(j+t)^{s_{pq}}}$ when $j<0$, where $(p,q)\in D(\lambda)$ is the $k$th component of $\theta_i$. Now, the weight $w^{\pm N}_{\bm{s};t}(l_i)$ of the path $l_i$ is defined to be the product of weights of all northeast edges along $l_i$. Here, we understand that $w^{\pm N}_{\bm{s};t}(l_i)=1$ if $\theta_i=\emptyset$. Moreover, we define the weight $w^{\pm N}_{\bm{s};t}(L)$ of $L\in \mathcal{E}^{\pm N}_{\lambda}$ by 
\[
    w^{\pm N}_{\bm{s};t}(L)=\prod^{s}_{i=1}w^{\pm N}_{\bm{s};t}(l_i).
\]

\begin{example}
    Let $\lambda=(4,3,3,2)$. Consider the $E$-rim decomposition $\Theta\in\rim^{\lambda}_E$ of $\lambda$ appearing in Example~\ref{ex:E_rim}. It is easy to see that $\tau_E(\Theta)=(123)\in \mathfrak{S}_4$ via the following $L=(l_1,l_2,l_3,l_4)\in \mathcal{E}^{\pm 3}_{(4,3,3,2)}$
    \begin{figure}[ht]
        \begin{center}
            \begin{tikzpicture} 
                \node at (0.5,1) {$1$};
                \node at (0.5,2) {$2$};
                \node at (0.5,3) {$3$};
                \node at (0.5,4) {$-3$};
                \node at (0.5,5) {$-2$};
                \node at (0.5,6) {$-1$};
                \node at (0.5,7) {$0$};
                \node at (1,0) {$C_4$};  
                \node at (2,0) {$C_3$};
                \node at (3,0) {$C_2$};
                \node at (4,0) {$C_1$};
                \node at (1,0.5) {$(1,1)$};  
                \node at (2,0.5) {$(2,1)$};
                \node at (3,0.5) {$(3,1)$};
                \node at (4,0.5) {$(4,1)$};
                \node at (2,7.5) {$(2,0)$};
                \node at (5,7.5) {$(5,0)$};
                \node at (7,7.5) {$(7,0)$}; 
                \node at (8,7.5) {$(8,0)$};
                \node at (2,8) {$D_4$};
                \node at (5,8) {$D_3$};
                \node at (7,8) {$D_2$}; 
                \node at (8,8) {$D_1$}; 
                \node at (1,1) {$\bullet$};
                \node at (1,2) {$\bullet$};
                \node at (1,3) {$\bullet$};
                \node at (1,4) {$\bullet$};   
                \node at (1,5) {$\bullet$};
                \node at (1,6) {$\bullet$};
                \node at (1,7) {$\bullet$};    
                \node at (2,1) {$\bullet$};
                \node at (2,2) {$\bullet$};
                \node at (2,3) {$\bullet$};  
                \node at (2,4) {$\bullet$};
                \node at (2,5) {$\bullet$};
                \node at (2,6) {$\bullet$};
                \node at (2,7) {$\bullet$};    
                \node at (3,1) {$\bullet$};
                \node at (3,2) {$\bullet$};
                \node at (3,3) {$\bullet$};  
                \node at (3,4) {$\bullet$};
                \node at (3,5) {$\bullet$};
                \node at (3,6) {$\bullet$};
                \node at (3,7) {$\bullet$};    
                \node at (4,1) {$\bullet$};
                \node at (4,2) {$\bullet$};
                \node at (4,3) {$\bullet$};  
                \node at (4,4) {$\bullet$};  
                \node at (4,5) {$\bullet$};
                \node at (4,6) {$\bullet$};
                \node at (4,7) {$\bullet$};    
                \node at (5,1) {$\bullet$};
                \node at (5,2) {$\bullet$};
                \node at (5,3) {$\bullet$};  
                \node at (5,4) {$\bullet$};  
                \node at (5,5) {$\bullet$};
                \node at (5,6) {$\bullet$};
                \node at (5,7) {$\bullet$};    
                \node at (6,1) {$\bullet$};
                \node at (6,2) {$\bullet$};
                \node at (6,3) {$\bullet$}; 
                \node at (6,4) {$\bullet$};  
                \node at (6,5) {$\bullet$};
                \node at (6,6) {$\bullet$};
                \node at (6,7) {$\bullet$};    
                \node at (7,1) {$\bullet$};
                \node at (7,2) {$\bullet$};
                \node at (7,3) {$\bullet$}; 
                \node at (7,4) {$\bullet$};  
                \node at (7,5) {$\bullet$};
                \node at (7,6) {$\bullet$};
                \node at (7,7) {$\bullet$};    
                \node at (8,1) {$\bullet$};
                \node at (8,2) {$\bullet$};
                \node at (8,3) {$\bullet$};  
                \node at (8,4) {$\bullet$};  
                \node at (8,5) {$\bullet$};
                \node at (8,6) {$\bullet$};
                \node at (8,7) {$\bullet$};    
                \draw (4,1) -- (5,2) -- (5,3) -- (5,4) -- (5,5) -- (6,6) -- (7,7);
                \draw[dotted] (3,1) -- (3,2) -- (3,3) -- (4,4) -- (4,5) -- (5,6) -- (5,7);
                \draw[loosely dashdotted] (2,1) -- (3,2) -- (4,3) -- (5,4) -- (6,5) -- (7,6) -- (8,7);      
                \draw[dashed] (1,1) -- (1,2) -- (1,3) -- (2,4) -- (2,5) -- (2,6) -- (2,7);
            \end{tikzpicture}
        \end{center}
        \caption{$L=(l_1,l_2,l_3,l_4)\in \mathcal{E}^{\pm 3}_{(4,3,3,2)}$}
    \end{figure}
    
    Let
    \(
        \ytableausetup{boxsize=18pt,aligntableaux=center}
        \bm{s}=
        \begin{ytableau}
        a & b & c & d \\
        e & f & g \\
        h & i & j \\
        k & l 
        \end{ytableau}
        \in T((4,3,3,2),\CC)
    \).
    Then, the weights of the $l_i$ are given by  
    \begin{align*}
        w^{\pm 3}_{\bm{s};t}(l_1)
        =&\frac{1}{1^a(-2+t)^e(-1+t)^h}, \quad
        w^{\pm 3}_{\bm{s};t}(l_2)
        =\frac{1}{3^b(-2+t)^f}, \\
        w^{\pm 3}_{\bm{s};t}(l_3)
        =&\frac{1}{1^c2^g3^j(-3+t)^i(-2+t)^l(-1+t)^k}, \quad
        w^{\pm 3}_{\bm{s};t}(l_4)
        =\frac{1}{3^d}.
    \end{align*}
    In particular, when 
    \(
        \ytableausetup{boxsize=18pt,aligntableaux=center}
        \bm{s}=
        \begin{ytableau}
        a_0 & a_1 & a_2 & a_3 \\
        a_{-1} & a_{0} & a_1 \\
        a_{-2} & a_{-1} & a_0 \\
        a_{-3} & a_{-2}
        \end{ytableau}
        \in T^{\diag}((4,3,3,2),\CC)
    \)
    , these are equal to 
    \begin{align*}
        &w^{\pm 3}_{\bm{s};t}(l_1)
        =\frac{1}{1^{a_0}(-2+t)^{a_{-1}}(-1+t)^{a_{-2}}}, \quad
        w^{\pm 3}_{\bm{s};t}(l_2)
        =\frac{1}{3^{a_1}(-2+t)^{a_0}}, \\
        &w^{\pm 3}_{\bm{s};t}(l_3)
        =\frac{1}{1^{a_2}2^{a_1}3^{a_0}(-3+t)^{a_{-1}}(-2+t)^{a_{-2}}(-1+t)^{a_{-3}}}, \quad
        w^{\pm 3}_{\bm{s};t}(l_4)
        =\frac{1}{3^{a_3}}.
    \end{align*}
    Notice that, in this case, from the definition of the weight, the tuple of indices of the exponent of the denominator of $w^{\pm 3}_{\bm{s};t}(l_i)$ along $l_i$ should be equal to $(a_{i-1},a_{i-2},a_{i-3},\ldots)$ for all $i$. 
\end{example}

We similarly give a proof of~\eqref{eq:jacobi_trudi_formula_Htype} by calculating the sum 
\[
    Y^{\pm N}_{\lambda}(\bm{s};t)
    =\sum_{L\in \mathcal{E}^{\pm N}_{\lambda}}\varepsilon_{\type(L)}w^{\pm N}_{\bm{s};t}(L)
    =\sum_{\sigma\in \mathfrak{S}^{\lambda}_E}
    \varepsilon_{\sigma}
    \sum_{L:C\to D^{\sigma}}
    w^{\pm N}_{\bm{s};t}(L).
\]

\begin{lemma} 
\label{lemma:E_step_1}
    For $\sigma\in \mathfrak{S}^{\lambda}_E$, let $\Theta^{\sigma}=(\theta^{\sigma}_1,\ldots,\theta^{\sigma}_s)\in \rim^{\lambda}_E$ be the $E$-rim decomposition such that $\tau_E(\Theta^{\sigma})=\sigma$. Then, we have  
    \[ 
        \sum_{L:C\to D^{\sigma}}w^{\pm N}_{\bm{s};t}(L) = \prod^{s}_{i=1}\zeta^{\pm N}\left(\theta^{\sigma}_i(\bm{s});t\right).
    \]
    Here, for $\Theta=(\theta_1,\ldots,\theta_s)\in\rim^{\lambda}_E$, $\theta_i(\bm{s})\in\CC^{|\theta_i|}$ is the tuple obtained by reading the entries of the restriction of $\bm{s}$ to $\theta_i$ from the top right to the bottom left.
\end{lemma}

From Lemma~\ref{lemma:E_step_1}, we have 
\begin{equation}
\label{eq:YN}
    Y^{\pm N}_{\lambda}(\bm{s};t)
    =\sum_{\sigma\in \mathfrak{S}^{\lambda}_E}
    \varepsilon_{\sigma}
    \prod^{s}_{i=1}\zeta^{\pm N}\left(\theta^{\sigma}_i(\bm{s});t\right).
\end{equation}
Define $\mathcal{E}^{\pm N}_{\lambda,0}$ similarly to $\mathcal{H}^{\pm N}_{\lambda,0}$ and also $ Y^{\pm N}_{\lambda,0}(\bm{s};t)$ and $Y^{\pm N}_{\lambda,1}(\bm{s};t)$. It holds that 
\[
    Y^{\pm N}_{\lambda,0}(\bm{s};t)
    =\sum_{L\in \mathcal{E}^{\pm N}_{\lambda,0}}w^{\pm N}_{\bm{s};t}(L)
    =\zeta^{\pm N}_{\lambda}(\bm{s};t).
\]
Hence, from~\eqref{eq:YN}, we reach the expression  
\begin{align}
\label{eq:Y_main_expression}
    \zeta^{\pm N}_{\lambda}(\bm{s};t) &= \sum_{\sigma\in \mathfrak{S}^{\lambda}_E} \varepsilon_{\sigma} \prod^{s}_{i=1} \zeta^{\pm N}\left(\theta^{\sigma}_i(\bm{s});t\right) - Y^{\pm N}_{\lambda,1}(\bm{s};t).
\end{align}

Now,~\eqref{eq:jacobi_trudi_formula_Htype} is obtained by taking the limit $N\to\infty$ of~\eqref{eq:Y_main_expression} under suitable assumptions on $\bm{s}$ described in Theorem~\ref{th:jacobi_trudi_formula} together with the following proposition.

\begin{prop}
\label{prop:jacobi_trudi_formula_Y}
    Assume that $\bm{s}=(s_{ij})\in T^{\diag}(\lambda,\CC)$. Write $a_{k}=s_{i,i+k}$ for $k\in\ZZ$.
    \begin{itemize}
        \item[$(1)$] 
        We have 
        \[
            Y^{\pm N}_{\lambda}(\bm{s};t) = \det \left[ \zeta^{\pm N}(a_{j-1},\ldots,a_{j-(\lambda'_i-i+j)};t) \right]_{1\leq i,j\leq s}.
        \]
        \item[$(2)$]
        It holds that
        \[
            Y^{\pm N}_{\lambda,1}(\bm{s};t)=0.
        \]
    \end{itemize}
\end{prop}

\begin{remark}
    When $\bm{s}\in T^{\diag}(\lambda,\CC)$,~\eqref{eq:Y_main_expression} can also be written in terms of the $E$-rim decomposition as follows:
    \begin{equation}
    \label{for:keyexpressionErim}
        \zeta^{\pm N}_{\lambda}(\bm{s};t) = \sum_{\Theta=(\theta_1,\ldots,\theta_s)\in \rim^{\lambda}_E} \varepsilon_E(\Theta) \zeta^{\pm N}\left(\theta_1(\bm{s});t\right)\cdots \zeta^{\pm N}\left(\theta_s(\bm{s});t\right),
    \end{equation}
    where $\varepsilon_E(\Theta)=\varepsilon_{\tau_E(\Theta)}$. Note that $\varepsilon_E(\Theta)=(-1)^{n-\#\{ i \mid \theta_i \neq \emptyset\}}$ when $\lambda=(n)$. 
\end{remark}

\section{Giambelli-type Formula}
\subsection{Expressions of unified Schur multiple zeta functions of hook type}
Let $p,q\in\ZZ_{\geq 0}$ and $\bm{s}=(s_{ij})\in T(\lambda,\CC)$. When $\lambda=(p+1,1^q)$, let $s_{ij}=z_{j-i}$, and write it as follows.
$$
\bm{s}=
{
\ytableausetup{boxsize=2.5em}
\begin{ytableau}
s_{11}&s_{12}& \cdots & s_{1,p+1}\\
s_{21}\\
\vdots\\
s_{q+1,1}
\end{ytableau}
}
=
{
\ytableausetup{boxsize=2.5em}
\begin{ytableau}
z_{0}&z_{1}& \cdots & z_{p}\\
z_{-1}\\
\vdots\\
z_{-q}
\end{ytableau}
}
\in W_{\lambda}.
$$
The unified Schur multiple zeta function of this type is called a
{\emph{unified Schur multiple zeta function of hook type}}, 
and we obtain the following expression for it.

\begin{thm}
For $\lambda=(p+1,1^q)$ and $\bm{s} \in W_{\lambda}$, we have
\begin{align}
	\zeta_{\lambda}^{\unifhat}(\bm{s};t) &= \sum_{j=0}^q (-1)^j \zeta^{\unifhat \star}(z_{-j},\ldots,z_{-1},z_0,z_1,\ldots,z_{p};t) \zeta^{\unifhat}(z_{-j-1},\ldots,z_{-q+1},z_{-q};t), \label{eq:schur_hook_type_expression} \\
	&= \sum_{j=0}^p (-1)^j \zeta^{\unifhat}(z_j,\ldots,z_2,z_1,z_0,z_{-1},\ldots,z_{-q};t) \zeta^{\unifhat \star}(z_{j+1},\ldots,z_{p-1},z_p;t).
\end{align}
\end{thm}

\begin{proof}
We adapt the proof of Theorem~3.1 in 
\cite{matsumoto2024schurgiambelli}. The difference is~\eqref{eq:schur_hook_type_expression_proof} below.
For the first assertion~\eqref{eq:schur_hook_type_expression}, we carry out the induction on $q$. In the case $q=0$, it is trivial that the assertion holds. 
Assume that the assertion holds for $q-1$, and consider the case of $q$:
\begin{multline}
\label{eq:schur_hook_type_expression_proof}
\zeta_{\lambda}^{\unifhat}(\bm{s};t)
= \sum_{\substack{m_{11}\preceq m_{12}\preceq\cdots\preceq m_{1,p+1}\\m_{11}\prec m_{21}\prec\cdots\prec m_{q+1,1}}}(m_{11}+\mathds{1}_{<0}(m_{11})t)^{-z_0}
(m_{12}+\mathds{1}_{<0}(m_{12})t)^{-z_1}\cdots \\
(m_{1,p+1}+\mathds{1}_{<0}(m_{1,p+1})t)^{-z_p} (m_{21}+\mathds{1}_{<0}(m_{21})t)^{-z_{-1}}\cdots(m_{q+1,1}+\mathds{1}_{<0}(m_{q+1,1})t)^{-z_{-q}}.
\end{multline}

Set $\lambda_{-}=(p+2,1^{q-1})$. The summation on the right-hand side can be divided into two parts:
\begin{equation*}
\sum_{\substack{m_{11}\preceq m_{12}\preceq\cdots\preceq m_{1,p+1}\\m_{21}\prec\cdots\prec m_{q+1,1}}}-\sum_{\substack{m_{11}\preceq m_{12}\preceq\cdots\preceq m_{1,p+1}\\m_{11}\succeq m_{21}\prec\cdots\prec m_{q+1,1}}}={\sum}_1-{\sum}_2.
\end{equation*}
Then, we have
$${\sum}_1=\zeta^{\unifhat \star}(z_0,\ldots,z_p;t)\zeta^{\unifhat}(z_{-1},\ldots,z_{-q};t),$$
while
$${\sum}_2=\sum_{\substack{m_{21}\preceq m_{11}\preceq m_{12}\preceq\cdots\preceq m_{1,p+1}\\m_{21}\prec\cdots\prec m_{q+1,1}}}=\zeta_{\lambda_{-}}^{\unifhat}(\bm{s}^{\flat};t),$$
where 
$\bm{s}^{\flat}=
{
\ytableausetup{boxsize=2.5em}
\begin{ytableau}
s_{21}&s_{11}&s_{12}& \cdots & s_{1,p+1}\\
s_{31}\\
\vdots\\
s_{q+1,1}
\end{ytableau}
}=
{
\ytableausetup{boxsize=2.5em}
\begin{ytableau}
z_{-1}&z_{0}&z_{1}& \cdots & z_{p}\\
z_{-2}\\
\vdots\\
z_{-q}
\end{ytableau}
}
$.

By the induction hypothesis, we obtain
 \begin{align*}
 {\sum}_2&=\sum_{i=0}^{q-1}(-1)^i\zeta^{\unifhat \star}(z_{-(i+1)},\ldots,z_{-2},z_{-1},z_0,\ldots,z_p;t)\zeta^{\unifhat}(z_{-(i+2)},\ldots,z_{-q};t)\\
 &=\sum_{j=1}^{q}(-1)^{j-1}\zeta^{\unifhat \star}(z_{-j},\ldots,z_{-2},z_{-1},z_0,\ldots,z_p;t)\zeta^{\unifhat}(z_{-(j+1)},\ldots,z_{-q};t).
 \end{align*}
This leads to
\begin{equation}
\begin{aligned}
{\sum}_1-{\sum}_2&=\zeta^{\unifhat \star}(z_0,\ldots,z_p;t)\zeta^{\unifhat}(z_{-1},\ldots,z_{-q};t)\\
&\quad -\sum_{j=1}^{q}(-1)^{j-1}\zeta^{\unifhat \star}(z_{-j},\ldots,z_{-2},z_{-1},z_0,\ldots,z_p;t)\zeta^{\unifhat}(z_{-(j+1)},\ldots,z_{-q};t)\\
&=\sum_{j=0}^{q}(-1)^{j}\zeta^{\unifhat \star}(z_{-j},\ldots,z_{-2},z_{-1},z_0,\ldots,z_p;t)\zeta^{\unifhat}(z_{-(j+1)},\ldots,z_{-q};t),
\end{aligned} 
\end{equation}
which completes the proof of~\eqref{eq:schur_hook_type_expression}. The second formula of the theorem can be similarly proved by induction on $p$.
\end{proof}

\subsection{Unified Schur multiple zeta functions of Hurwitz type}

In this section, we give the Giambelli-type formula for unified Schur multiple zeta functions, by a method quite similar to that developed in \cite{matsumoto2025schurgiambelli} and \cite{matsumoto2025schur}.

\begin{definition}
For $\bm{s}=(s_{ij})\in W_{\lambda},\;\bm{t}=(t_{ij})\in T({\lambda},\RR_{<1})$, we define
$$\zeta_{\lambda}^{\unifhat}(\bm{s}|\bm{t}) := \lim_{N\to\infty}\sum_{M=(m_{ij})\in\SSYT^{\pm N}({\lambda})}\prod_{(i,j)\in D(\lambda)}\frac{1}{(m_{ij}+\mathds{1}_{<0}(m_{ij})t_{ij})^{s_{ij}}}.$$
\end{definition}

The absolute convergence of the series is similar to the one defining $\zeta_{\lambda}^{\unifhat}(\bm{s};t)$. Note that when all the entries of $\bm{t}$ are equal, we clearly have $\zeta_{\lambda}^{\unifhat}(\bm{s}|\bm{t}) = \zeta_{\lambda}^{\unifhat}(\bm{s};t)$.

For $\bm{s}\in T(\lambda, \CC)$ and  $X=(x_{ij})\in T(\lambda,\RR)$, let
$$
Z_\lambda (\bm{s},X):=\prod_{(i,j)\in D(\lambda)}\frac{1}{x_{ij}^{s_{ij}}}.
$$
Therefore,
\begin{align*}
\zeta_{\lambda}^{\unifhat}(\bm{s}|\bm{t})&=\lim_{N\to\infty}\sum_{M\in \SSYT^{\pm N}(\lambda)}\prod_{(i,j)\in D(\lambda)}\frac{1}{(m_{ij}+\mathds{1}_{<0}(m_{ij})t_{ij})^{s_{ij}}}\\
&=\lim_{N\to\infty}\sum_{M\in \SSYT^{\pm N}(\lambda)}Z_{\lambda}(\bm{s},M+\mathds{1}_{<0}(M)\bm{t}).
\end{align*}

We now introduce another notation for $\lambda$ called Frobenius notation.
For a partition $\lambda=(\lambda_1,\lambda_2,\ldots),$ let $\lambda^{\prime}=(\lambda_1^{\prime},\lambda_2^{\prime},\ldots)$ be its conjugate, and we define two sequences of indices $p_1,\ldots,p_n$ and $q_1,\ldots,q_n$ by $p_i=\lambda_i-i$ and $q_i=\lambda_i^{\prime}-i$ for $1\le i\le n$ where $n$ is the number of entries in the main diagonal of the Young diagram of $\lambda$. We sometimes write $\lambda=(p_1,\ldots,p_n|q_1,\ldots,q_n)$, the Frobenius notation of $\lambda$. Let $\alpha(\lambda)\subset\lambda$ be the set of all cells which are to the right of the northeast boundary of the main diagonal in the Young diagram of shape $\lambda$ and $\beta(\lambda)\subset\lambda$ be the set of all cells which are to the left of the southwest boundary of the main diagonal in the Young diagram of shape $\lambda$. We denote them by $\alpha(\lambda)=(p_1,\ldots,p_n),\beta(\lambda)=(q_1,\ldots,q_n)$. Moreover, let $D^M(\lambda)$ be the set of all cells which are on the main diagonal in the Young diagram of shape $\lambda$.

For $\sigma\in\mathfrak{S}_n$ and a Young tableau $X=(x_{ij}) \in T(\lambda,\llbracket -N,N \rrbracket \setminus \{ 0 \}\})$, we call $X$ a $\sigma$-tableau\footnote{The integer $N$ will be implicit in this section.} if it satisfies the following conditions:
\begin{enumerate}[label=(\roman*)]
	\item $x_{i,j}\preceq x_{i,j+1}$ along the rows in $\alpha(\lambda)$.
	\item $x_{i,j}\prec x_{i+1,j}$ down the columns in $\beta(\lambda)\cup D^M(\lambda)$.
	\item $x_{\sigma(i),\sigma(i)}\preceq x_{i,i+1}$ for $i=1,2,\ldots,n$ whenever $i+1\le\lambda_i$.
\end{enumerate}

For $\bm{s}=(s_{ij})_{(i,j) \in D(\lambda)} \in T^{\diag}(\lambda,\CC)$, define $\bm{s}_{ij}^F=
{
	\ytableausetup{boxsize=2em}
	\begin{ytableau}
		z_0 & z_1 & z_2 & \cdots  & z_{p_i}\\
		z_{-1}\\
		\vdots\\
		z_{-q_j}
	\end{ytableau}
}$.

\begin{lemma}\label{lem:expression_Z_hook_type}
Let $X=(x_{ij})$ be a $\sigma$-tableau of shape $\lambda=(p_1,\ldots, p_n|q_1,\ldots, q_n).$ Then,
\begin{align*}
X_k(\sigma)=
  \begin{array}{|c|c|c|c|}
 \hline
 x_{\sigma(k),\sigma(k)} & x_{k,k+1} & \cdots & x_{k,k+p_k}\\
 \hline
 x_{\sigma(k)+1,\sigma(k)}\\
 \cline{1-1}
 \vdots\\
 \cline{1-1}
 x_{\sigma(k)+q_{\sigma(k)},\sigma(k)}\\
 \cline{1-1}
  \end{array}
  \; 
\end{align*}
is in $\SSYT^{\pm N}({p_k|q_{\sigma(k)}})$ for $1\le k\le n$, and for $\bm{s}\in W^{\diag}_{\lambda}$ and $\bm{t}\in T^{\diag}(\lambda,\RR_{<1})$,
\begin{equation}
Z_{\lambda}(\bm{s},X+\mathds{1}_{<0}(X)\bm{t})=\prod_{k=1}^n Z(\bm{s}_{k,\sigma(k)}^F, X_k(\sigma)+\mathds{1}_{<0}(X_k(\sigma))\bm{t}_{k,\sigma(k)}^F).
\end{equation}
\end{lemma} 

\begin{proof}
By the condition (i) 
of $\sigma$-tableau, $x_{k, k+1}\preceq\cdots\preceq x_{k,k+p_k}$
, and by the condition (ii), $x_{\sigma(k),\sigma(k)}\prec x_{\sigma(k)+1,\sigma(k)}\prec\cdots\prec x_{\sigma(k)+q_{\sigma(k)},\sigma(k)}$. Furthermore, by the condition (iii) of $\sigma$-tableau, $x_{\sigma(k),\sigma(k)}\preceq x_{k,k+1}$. So $X_k(\sigma)$ belongs to $\SSYT^{\pm N}$.

By the definition of $Z_{\lambda}(\bm{s},X+\mathds{1}_{<0}(X)\bm{t})$, we have
\begin{equation*}
\begin{aligned}
&Z_{\lambda}(\bm{s},X+\mathds{1}_{<0}(X)\bm{t})=\prod_{(i,j)\in D(\lambda)}\frac{1}{(x_{ij}+\mathds{1}_{<0}(x_{ij})t_{ij})^{s_{ij}}}\\
&=\prod_{(i,j)\in \alpha(\lambda)}\frac{1}{(x_{ij}+\mathds{1}_{<0}(x_{ij})t_{ij})^{s_{ij}}}\prod_{(i,j)\in \beta(\lambda)\cup D^M(\lambda)}\frac{1}{(x_{ij}+\mathds{1}_{<0}(x_{ij})t_{ij})^{s_{ij}}}\\
&=\prod_{(i,j)\in D(p_1)}\frac{1}{(x_{ij}+\mathds{1}_{<0}(x_{ij})t_{ij})^{s_{ij}}}\cdots\prod_{(i,j)\in D(p_n)}\frac{1}{(x_{ij}+\mathds{1}_{<0}(x_{ij})t_{ij})^{s_{ij}}}\\
&\quad\times\prod_{(i,j)\in D(q_{\sigma(1)})\cup (\sigma(1),\sigma(1))}\frac{1}{(x_{ij}+\mathds{1}_{<0}(x_{ij})t_{ij})^{s_{ij}}}\cdots\\
&\qquad\cdots\prod_{(i,j)\in D(q_{\sigma(n)})\cup(\sigma(n),\sigma(n))}\frac{1}{(x_{ij}+\mathds{1}_{<0}(x_{ij})t_{ij})^{s_{ij}}}\\
&=\left(\prod_{(i,j)\in D(p_1)}\frac{1}{(x_{ij}+\mathds{1}_{<0}(x_{ij})t_{ij})^{s_{ij}}}\prod_{(i,j)\in D(q_{\sigma(1)})\cup (\sigma(1),\sigma(1))}\frac{1}{(x_{ij}+\mathds{1}_{<0}(x_{ij})t_{ij})^{s_{ij}}}\right)\\
&\quad\cdots\left(\prod_{(i,j)\in D(p_n)}\frac{1}{(x_{ij}+\mathds{1}_{<0}(x_{ij})t_{ij})^{s_{ij}}}\prod_{(i,j)\in D(q_{\sigma(n)})\cup (\sigma(n),\sigma(n))}\frac{1}{(x_{ij}+\mathds{1}_{<0}(x_{ij})t_{ij})^{s_{ij}}}\right)\\
&=\prod_{k=1}^nZ(\bm{s}_{k,\sigma(k)}^F, X_k(\sigma)+\mathds{1}_{<0}(X_k(\sigma))\bm{t}_{k,\sigma(k)}^F).
\end{aligned}
\end{equation*}
\end{proof}
Conversely, for any given $X_k(\sigma)\in\SSYT^{\pm N}(p_k|q_{\sigma(k)})(1\le k\le n)$, we may construct a $\sigma$-tableau $X$ by combining those $X_k(\sigma)$ as in the above proof.

Let
$$
\mathbb{X}:=\{(\sigma,X) \mid \sigma\in\mathfrak{S}_n, X \in \SSYT^{\pm N}(\lambda):\sigma -\text{tableau of shape }\lambda\}.
$$
We see that $\mathbb{X}=\mathbb{X}^{\dag}\cup\mathbb{X}^{\dag\dag}$, where $\mathbb{X}^{\dag}:=\{(\id
,X) \mid X\in\SSYT^{\pm N}(\lambda)\}$ where $\id$ is the identity permutation, and $\mathbb{X}^{\dag\dag}=\mathbb{X}\backslash\mathbb{X}^{\dag}$.
If  $X\in\SSYT^{\pm N}(\lambda)$, then $X$ is obviously an $\id$-tableau. However, $X$ is never a $\sigma$-tableau for any $\sigma\neq \id$.
In fact, if $X$ is a $\sigma$-tableau, then $x_{\sigma(i),\sigma(i)}\preceq x_{i,i+1}$ for any $i$. If $\sigma\neq \id$,
then there exists an $i$ for which $j=\sigma(i)>i$, so $x_{j,j}\preceq x_{i,i+1}$. However, since $X$ is $\SSYT^{\pm N}$, it should be
$x_{i,i+1}\preceq x_{j-1,j}\prec x_{j,j}$, which is a contradiction.
\begin{lemma}
\label{lemma:relation_z_x_xprime}
Let $\bm{t}\in T^{\diag}(\lambda,\RR_{<1})$ and $\bm{s} \in W_{\lambda}^{\diag}$. For any $(\sigma, X)\in\mathbb{X}^{\dag\dag}$, there exists another element $(\sigma^{\prime}, X^{\prime})\in\mathbb{X}^{\dag\dag}$ such that
\begin{equation}\label{eq:formula_signed_x_xprime}
\sgn(\sigma)Z_{\lambda}(\bm{s},X+\mathds{1}_{<0}(X)\bm{t})+\sgn(\sigma^{\prime})Z_{\lambda}(\bm{s},X^{\prime}+\mathds{1}_{<0}(X')\bm{t})=0.
\end{equation}
\end{lemma}
\begin{proof}
Let $(\sigma,X)\in \mathbb{X}^{\dag\dag}$. Then $X\notin\SSYT^{\pm N}(\lambda)$. The situation can be divided into the following three cases.
\begin{enumerate}[label=(\Roman*)]
	\item $X$ does not satisfy the $\SSYT^{\pm N}$ condition in the region $\alpha(\lambda)$;
	\item $X$ satisfies the $\SSYT^{\pm N}$ condition in $\alpha(\lambda)$, but does not satisfy it in the region $\beta(\lambda)\cup D^M(\lambda)$;
	\item $X$ satisfies $\SSYT^{\pm N}$ condition both in $\alpha(\lambda)$ and in $\beta(\lambda)\cup D^M(\lambda)$, but does not satisfy it on the north-east boundary of the main diagonal.
\end{enumerate}
We treat these three cases separately.\\
{\it{Case(I).}}
Assume that the SSYT$^{\pm N}$ condition is violated in $\alpha(\lambda)$, between the $i$-th row and the $(i+1)$-th row. We can find a $k(\ge i+2)$ for which 
$x_{i,k}\succeq x_{i+1,k}$. 
Let $k$ be the largest of such integers. Let $\sigma^{\prime}=\sigma\circ (i,i+1)$, and define $X'=(x'_{ij})$ by interchanging the entries $x_{i,h}$ and $x_{i+1,h+1}$ of $X\;(i+1\le h\le k-1)$, while all other entries remain the same as $X$.

We show that $X'$ is a $\sigma '$-tableau. 
In fact,
$x_{i,j}\preceq x_{i,j+1}$ along the $i$-th row of $\alpha(\lambda)$ in $X'$.
Because $x_{i+1,k}\preceq x_{i,k}$ and
$x_{i,k-1} \preceq x_{i,k} \preceq x_{i,k+1}\prec x_{i+1,k+1}$ in the $(i+1)$-th row of 
$X'$ (where the last inequality comes from the maximality of $k$), the condition (i) is valid.
The condition (ii) is obviously valid because $X'$ in $\beta(\lambda)\cup D^M(\lambda)$ is exactly the same as $X$.

Next we check the condition (iii). Since $\sigma'=\sigma\circ(i,i+1)$, we find
\begin{align*}
&x'_{\sigma'(i),\sigma'(i)}=x'_{\sigma(i+1),\sigma(i+1)}=x_{\sigma'(i),\sigma'(i)}\preceq x_{i+1,i+2}=x'_{i,i+1},\\
&x'_{\sigma'(i+1),\sigma'(i+1)}=x'_{\sigma(i),\sigma(i)}=x_{\sigma(i),\sigma(i)}\preceq x_{i,i+1}=x'_{i+1,i+2}
\end{align*}
and for any $l\neq i,i+1$,
$$
x'_{\sigma'(i),\sigma'(i)}=x'_{\sigma(i),\sigma(i)}=x_{\sigma(i),\sigma(i)}\preceq x_{l,l+1}=x'_{l,l+1},
$$
hence (iii) is valid. Therefore $X'$ is a $\sigma'$-tableau.
Since $\bm{s}\in W_{\lambda}^{\diag}$, we see that
\begin{equation}\label{eq:equality_z_x_xprime}
Z_{\lambda}(\bm{s},X)=Z_{\lambda}(\bm{s},X').
\end{equation}
We now consider the case of $X+\mathds{1}_{<0}(X)\bm{t}$ and $X'+\mathds{1}_{<0}(X')\bm{t}$ instead of $X$ and $X'$ in~\eqref{eq:equality_z_x_xprime}. Since $\bm{t}\in T^{\diag}(\lambda,\RR_{<1})$, the same element is added to the section changed between $X$ and $X'$, respectively. Therefore, we have 
\[
	Z_{\lambda}(\bm{s},X+\mathds{1}_{<0}(X)\bm{t})=Z_{\lambda}(\bm{s},X'+\mathds{1}_{<0}(X')\bm{t}).
\]
Noting $\sgn(\sigma')=-\sgn(\sigma)$, we obtain~\eqref{eq:formula_signed_x_xprime} in this case.

{\it{Case(II).}}
Assume that the SSYT$^{\pm N}$ condition is valid in $\alpha(\lambda)$, but is violated in $\beta(\lambda)\cup D^M(\lambda)$, between the $j$-th column and the $(j+1)$-th column. We can find a $k(\ge j+1)$ for which $x_{k,j}\succ x_{k,j+1}$. We may assume that $k$ is the largest such $k$.
Let $\sigma'=(j,j+1)\circ\sigma$, and this time define $X'=(x'_{ij})$ by changing the places of the elements $x_{h,j}$ and $x_{h+1,j+1}$ of $X\;(j\le h\le k-1)$, while all the other elements remain the same as $X$.

Again, we can show that $X'$ is a $\sigma'$-tableau. The condition (i) is obvious, while the condition (ii) follows from $x_{k,j+1}\prec x_{k,j}$ and $x_{k-1,j}\prec x_{k,j}\prec x_{k+1,j}\preceq x_{k+1,j+1}$. Let us check the condition (iii). For any $m$, we have 
\begin{equation*}
    \sigma'(m)=
    \begin{cases}
        \sigma(m) &\text{if } \sigma(m)\neq j,j+1,\\
    	j+1 &\text{if } \sigma(m)=j,\\
    	j &\text{if } \sigma(m)=j+1.
	\end{cases}
\end{equation*}
Therefore we find
\begin{equation*}
    x'_{\sigma'(k),\sigma'(k)}=
	\begin{cases}
    	x'_{\sigma(k),\sigma(k)}=x_{\sigma(k),\sigma(k)}\preceq x_{k,k+1}=x'_{k,k+1} &\text{if } \sigma(k)\neq j,j+1,\\
    	x'_{j+1,j+1}=x_{j,j}=x_{\sigma(k),\sigma(k)}\preceq x_{k,k+1}=x'_{k,k+1} &\text{if } \sigma(k)=j,\\
    	x'_{j,j}=x_{j+1,j+1}=x_{\sigma(k),\sigma(k)}\preceq x_{k,k+1}=x'_{k,k+1} &\text{if } \sigma(k)=j+1.
	\end{cases}
\end{equation*}
which implies the condition (iii). Therefore $X'$ is $\sigma'$-tableau. Since $\bm{s}\in W_{\lambda}^{\diag}$, $\bm{t}\in T^{\diag}(\lambda,\RR_{<1})$, and $\sgn(\sigma')=-\sgn(\sigma)$,~\eqref{eq:formula_signed_x_xprime} follows in this case.

{\it{Case(III).}}
Assume that the SSYT$^{\pm N}$ condition is valid in $\alpha(\lambda)$, also in $\beta(\lambda)\cup D^M(\lambda)$, but is violated on the north-east boundary of the main diagonal. Let $k$ be the least such that violation happens at $x_{k,k+1}$; that is, the first violation is either (A) $x_{k,k}\succ x_{k,k+1}$ or (B) $x_{k+1,k+1}\preceq x_{k,k+1}$. 

First we claim that $\sigma(h)=h$ for all $h\le k-1$ (when $k\ge 2$). Put $r=\sigma(1)$. Then $x_{r,r}=x_{\sigma(1),\sigma(1)}\preceq x_{1,2}$. But if $r\ge 2$, then $x_{r,r}\succeq x_{2,2}\succ x_{1,2}$ (because of the SSYT$^{\pm N}$ condition in $\beta(\lambda)\cup D^M(\lambda)$ and the minimality of $k$), and hence $r=1$. When $k=2$, we are done. Let $k\ge 3$, and put $s=\sigma(2)$. 
We have $x_{s,s}=x_{\sigma(2),\sigma(2)}\preceq x_{2,3}$, while if $s\ge 3$ then $x_{s,s}\succeq x_{3,3}\succ x_{2,3}$, so $s=2$. Repeating this procedure, we obtain the claim.

From the above claim we see that $\sigma(k)\ge k$. Therefore $x_{k,k}\preceq x_{\sigma(k),\sigma(k)}\preceq x_{k,k+1}$, so the case (A) does not happen.

Consider the case (B). Define $i,j$ by $\sigma(i)=k$ and $\sigma(j)=k+1$. By the above claim we have $i,j\ge k$. Put $\sigma'=\sigma\circ (i,j)$. Then
\begin{equation}\label{eq:ineq_caseIII_B}
x_{k+1,k+1}\preceq x_{k,k+1}\preceq x_{i,i+1}
\end{equation}
because of (B) and the SSYT$^{\pm N}$ condition in $\alpha(\lambda)$. Then
\begin{equation}\label{eq:ineq_caseIII_B_bis}
x_{k,k}\preceq x_{k+1,k+1}=x_{\sigma(j),\sigma(j)}\preceq x_{j,j+1}
\end{equation}
because of (B) and the SSYT$^{\pm N}$ condition in $\beta(\lambda)\cup D^M(\lambda)$.

In this case we choose $X=X'$, and prove that $X$ is also a $\sigma'$-tableau. The conditions (i) and (ii) are obvious. By~\eqref{eq:ineq_caseIII_B} and~\eqref{eq:ineq_caseIII_B_bis} we have
\begin{align*}
&x_{\sigma'(i),\sigma'(i)}=x_{\sigma(j),\sigma(j)}=x_{k+1,k+1}\preceq x_{i,i+1},\\
&x_{\sigma'(j),\sigma'(j)}=x_{\sigma(i),\sigma(i)}=x_{k,k}\preceq x_{j,j+1},
\end{align*}
and further,
$$
x_{\sigma'(l),\sigma'(l)}=x_{\sigma(l),\sigma(l)}\preceq x_{l,l+1}
$$
for any $l\neq i,j$. Therefore the condition (iii) is valid. So $X$ is a $\sigma'$-tableau. Since $\bm{s}\in W_{\lambda}^{\diag}$, $\bm{t}\in T^{\diag}(\lambda,\RR_{<1})$, and $\sgn(\sigma')=-\sgn(\sigma)$, we obtain~\eqref{eq:formula_signed_x_xprime}.
\end{proof}

The following follows immediately from Lemma~\ref{lemma:relation_z_x_xprime}.

\begin{cor}\label{cor:sum_sigma_tableau_dagger}
We have
\begin{align}
\sum_{(\sigma,X)\in\mathbb{X}}\sgn(\sigma)Z_{\lambda}(\bm{s},X+\mathds{1}_{<0}(X)\bm{t})&=\sum_{(\id,X)\in\mathbb{X}^{\dag}}\sgn(\id)Z_{\lambda}(\bm{s},X+\mathds{1}_{<0}(X)\bm{t})\\
&=\sum_{X\in\SSYT^{\pm N}(\lambda)}Z_{\lambda}(\bm{s},X+\mathds{1}_{<0}(X)\bm{t}).\notag
\end{align}
\end{cor}

\subsection{Giambelli formula for unified Schur--Hurwitz multiple zeta functions}

\begin{thm}
Let $\lambda$ be a partition whose Frobenius notation is $\lambda =(p_1,\ldots,p_n|q_1,\ldots,q_n)$. Assume $\bm{s}\in W^{\diag}_{\lambda}, \bm{t}\in T^{\diag}(\lambda ,\RR_{<1})$. Then we have
\begin{equation}\label{eq:giambelli_formula}
\zeta_{\lambda}^{\unifhat}(\bm{s}|\bm{t}) = \det(\zeta^{\unifhat}_{i,j}(\bm{s}_{ij}^F|\bm{t}_{ij}^F))_{1\le i,j\le n}
\end{equation}
where $\zeta^{\unifhat}_{i,j}(\bm{s}_{ij}^F|\bm{t}_{ij}^F) = \zeta^{\unifhat}_{(p_i+1, 1^{q_j})}(\bm{s}_{ij}^F|\bm{t}_{ij}^F)$.
\end{thm}
\begin{proof}
From Lemma~\ref{lem:expression_Z_hook_type}, the right-hand side of~\eqref{eq:giambelli_formula} is
\begin{align*}
&\det(\zeta^{\pm N}_{i,j}(\bm{s}_{ij}^F|\bm{t}_{ij}^F))_{1\le i,j\le n}=\sum_{\sigma\in\mathfrak{S}_n}\sgn(\sigma)\prod_{k=1}^n\zeta_{k,\sigma(k)}^{\pm N}(\bm{s}_{k,\sigma(k)}^F|\bm{t}_{k,\sigma(k)}^F)\\
&=\sum_{\sigma\in\mathfrak{S}_n}\sgn(\sigma)\prod_{k=1}^n\sum_{X_k(\sigma)\in \SSYT^{\pm N}(p_k|q_{\sigma(k)})}Z_{(p_k|q_{\sigma(k)})}(\bm{s}_{k,\sigma(k)}^F, X_k(\sigma)+\mathds{1}_{<0}(X_k(\sigma))\bm{t}_{k,\sigma(k)}^F)\notag\\
&=\sum_{\sigma\in\mathfrak{S}_n}\sgn(\sigma)\sum_{X:\sigma -\text{tableau}}Z_{\lambda}(\bm{s},X+\mathds{1}_{<0}(X)\bm{t})\notag.
\end{align*}
From the previous equality and Corollary~\ref{cor:sum_sigma_tableau_dagger}, we find
\begin{align}
\det(\zeta^{\pm N}_{i,j}(\bm{s}_{ij}^F|\bm{t}_{ij}^F))_{1\le i,j\le n}&=\sum_{(\sigma,X)\in \mathbb{X}}\sgn(\sigma)Z_{\lambda}(\bm{s},X+\mathds{1}_{<0}(X)\bm{t})\notag\\
&=\sum_{(\id,X)\in \mathbb{X}^{\dag}}\sgn(\id)Z_{\lambda}(\bm{s},X+\mathds{1}_{<0}(X)\bm{t})\notag\\
&=\sum_{X\in\SSYT^{\pm N}(\lambda)}Z_{\lambda}(\bm{s},X+\mathds{1}_{<0}(X)\bm{t}).\notag
\end{align}
Therefore,
$$
\lim_{N\to\infty}\det(\zeta^{\pm N}_{i,j}(\bm{s}_{ij}^F|\bm{t}_{ij}^F))_{1\le i,j\le n}=\lim_{N\to\infty}\sum_{X\in\SSYT^{\pm N}(\lambda)}Z_{\lambda}(\bm{s},X+\mathds{1}_{<0}(X)\bm{t})=\zeta_{\lambda}^{\unifhat}(\bm{s}|\bm{t}).
$$
\end{proof}

\section*{Acknowledgments}

Both authors wish to thank Maki Nakasuji for helpful discussions and comments. The second-named author is also grateful to Maki Nakasuji for a fruitful visit to Sophia University in August 2026.

\section*{AI use}

We acknowledge the use of ChatGPT 6 Astra to check for typos, improve the English throughout the text, and check some of the proofs in this paper. None of the ideas developed in this text are AI-generated. We take full responsibility for the content of the paper. This declaration is not an endorsement of OpenAI.

\nocite{matsumoto2025schur}

\bibliography{./Bibliography.bib}
\bibliographystyle{alpha}
\end{document}